\documentclass[a4paper,12pt]{article}
\usepackage{amsmath}
\usepackage{graphicx}
\usepackage{amssymb,geometry,latexsym}
\usepackage{citesort}
\usepackage{xypic}

\author{ Timothy Porter}

\title{$\mathcal{S}$-categories, $\mathcal{S}$-groupoids, Segal categories and quasicategories}

\newtheorem{theorem}{Theorem} 
\newtheorem{proposition}{Proposition} 
\newtheorem{corollary}{Corollary} 
\newtheorem{lemma}{Lemma} 
  
\begin{document}

\maketitle
  
The notes were prepared for a series of talks that I gave in Hagen in late June and early July 2003, and, with some changes, in the University of La Lagu\~
{n}a, the Canary Islands, in September, 2003. They assume the audience knows some abstract homotopy theory and as Heiner Kamps was in the audience in Hagen, it is safe to assume that the notes assume a reasonable knowledge of our book,  \cite{K&P}, or any equivalent text if one can be found!

What do the notes set out to do?

\medskip

\section*{``Aims and Objectives!'' or should it be ``Learning Outcomes''?}

\begin{itemize}
\item To revisit some oldish material on abstract homotopy and simplicially enriched categories, that seems to be being used in today's resurgence of interest in the area and to try to view it in a new light, or perhaps from  new directions;
\item To introduce Segal categories and various other tools used by the Nice-Toulouse group of abstract homotopy theorists and link them into some of the older ideas;
\item To introduce Joyal's quasicategories, (previously called weak Kan complexes but I agree with Andr\'e that his nomenclature is better so will adopt it) and show how that theory links in with some old ideas of Boardman and Vogt, Dwyer and Kan, and Cordier and myself;
\item To ask lots of questions of myself and of the reader.
\end{itemize} 
The notes include some material from the `Cubo' article, \cite{cubo}, which was itself based on notes for a course at the \emph{Corso  estivo  Categorie  e  Topologia} in 1991, but the overlap has been kept as small as is feasible as the purpose and the audience of the two sets of notes are different and the abstract homotopy theory has `moved on', in part, to try the new methods out on those same `old' problems and to attack new ones as well.

As usual when you try to specify `learning outcomes' you end up asking who has done the learning, the audience? Perhaps.  The lecturer, most certainly!

\vspace{1cm}

\textbf{Acknowledgements}

I would like to thank Heiner Kamps and his colleagues at the Fern Univerist\"{a}t for the invitation to give the talks of which these notes are a summary  and to the Fern Univerist\"{a}t for the support that made the visit possible, to Jos\'e Manuel Garc\'{\i}a-Calcines, Josu\'e Remedios and their colleagues and for the Departamento de Mathematica Fundamental in the Universidad de La Laguna, Tenerife,  simlarly and also to Carlos Simpson, Bertrand Toen, Andr\'e Joyal, Clemens Berger, Andr\'e Hirschowitz and others at the Nice meeting in May 2003, since that is where bits of ideas that I had gleaned over a longish period of time fitted together so that I think I begin to understand the way that a lot of things interlock in this area better than I did before! 

These notes have also benefitted from comments by Jim Stasheff and some of his colleagues on an earlier version.

\tableofcontents

\section{$\mathcal{S}$-categories}
\subsection{Categories with simplicial `hom-sets'}

We assume we have a category $\mathcal{A}$ whose objects will be denoted by lower case letter, $x$,$y$,$z$, \ldots , at least in the generic case, and for each pair of such objects, $(x,y)$, a simplicial set $\mathcal{A}(x,y)$ is given; for each triple $x, y, z$ of objects of $\mathcal{A}$, we have a simplicial map, called \emph{composition}
$$\mathcal{A}(x,y)\times \mathcal{A}(y,z)\longrightarrow \mathcal{A}(x,z);$$
and for each object $x$ a map 
$$\Delta[0] \to \mathcal{A}(x,x)$$
that `names' or `picks out' the `identity arrow' in the set of 0-simplices of $\mathcal{A}(x,x)$. This data is to satisfy the obvious axioms, associativity and  identity,   suitably adapted to this situation.  Such a set up will be called a \emph{simplicially enriched category} or more simply \emph{an $\mathcal{S}$-category}.  Enriched category theory is a well established branch of category theory.  It has many useful tools and not all of them have yet been exploited for the particular case of $\mathcal{S}$-categories and its applications in homotopy theory.

\medskip
{\sc Warning:} Some authors use the term simplicial category for what we have termed a simplicially enriched category.  There is a close link with  the notion of simplicial category that is consistent with usage in simplicial theory \emph{per se}, since any simplicially enriched category can be thought of as a simplicial object in the `category of categories', but a simplicially enriched category is not just a simplicial object in the `category of categories' and not all such simplicial objects correspond to such enriched categories. That being said that usage need not cause problems provided the reader is aware of the usage in the paper to which reference is being made.

\medskip

\textbf{Examples}

(i) $\mathcal{S}$, the category of simplicial sets:\\
here
$$\mathcal{S}(K,L)_n:= S(\Delta[n] \times K,L);$$
Composition : for $f \in \mathcal{S}(K,L)_n$, $ g \in \mathcal{S}(L,M)_n$, so $f : \Delta[n] \times K \to L$,  $g : \Delta[n] \times L \to M$,
$$g\circ f := (\Delta[n] \times K \stackrel{diag \times K}{\longrightarrow} \Delta[n] \times\Delta[n] \times K\stackrel{\Delta[n] \times f}{\longrightarrow }\Delta[n] \times L \stackrel{g}{\to} M);$$
Identity : $id_K : \Delta[0] \times K \stackrel{\cong}{\to} K$,

\medskip
(ii) $\mathcal{T}op$, `the' category of spaces (of course, there are numerous variants but you can almost pick whichever one you like as long as the constructions work): 
$$\mathcal{T}op(X,Y)_n := Top(\Delta^n \times X, Y)$$
Composition and identities are defined analogously to in (i).

\medskip

(iii) For each  $X$, $Y \in Cat$, the category of small categories, then
we similarly get $\mathcal{C}at(X,Y)$,
$$\mathcal{C}at(X,Y)_n = {Cat}([n] \times X, Y).$$
We leave the other structure up to the reader.

(iv) $\mathcal{C}rs$, the category of crossed complexes: see \cite{K&P} for background and other references, and Tonks, \cite{andythesis}
for a more detailed treatment of the simplicially enriched  category structure;
$$\mathcal{C}rs(A,B) := Crs(\pi(n)\otimes C, D)$$
Composition has to be defined using an approximation to the identity, again see \cite{andythesis}.

\medskip

(v) $\mathcal{C}h^+_K$, the category of positive chain complexes of modules over a commutative ring $K$. (Details are left to the reader, or follow from the Dold-Kan theorem and example (vi) below.)

\medskip

(vi) $\mathcal{S}(Mod_K)$, the category of simplicial $K$-modules.  The structure uses tensor product with the free simplicial $K$-module on $\Delta[n]$ to define the `hom' and the composition, so is very much like (i).

\medskip

In general any category of simplicial objects in a `nice enough' category has a simplicial enrichment, although the general argument that gives the construction does not always make the structure as transparent as it might be.

There is an evident notion of $\mathcal{S}$-enriched functor, so we get a category of `small' $\mathcal{S}$-categories, denoted $\mathcal{S}\!-\!Cat$.  Of course, none of the above examples are `small'.  (With regard to `smallness', although sometimes a smallness condition is essential, one can often ignore questions of smallness and, for instance, consider simplicial `sets' where actually the collections of simplices are not truly `sets' (depending on your choice of methods for handling such foundational questions).)

\subsection{From simplicial resolutions to $\mathcal{S}$-cats.}
  The forgetful functor $U : Cat \to DGrph_0$ has a left adjoint, $F$.  Here $DGrph_0$ denotes the category of directed graphs with `identity  loops', so $U$ forgets just the composition within each small category but remembers that certain loops are special `identity loops'.  The free category functor here takes, between any two objects, all strings of composable \emph{non-identity} arrows that start at the first object and end at the second. One can think of $F$ identifying the old identity arrow at an object $x$ with the empty string at $x$. 

This adjoint pair gives a comonad on $Cat$ in the usual way, and hence a functorial simplicial resolution, which we will denote $S(\mathbb{A})\to \mathbb{A}$.  In more detail, we write $T = FU$ for the functor part of the comonad, the unit of the adjunction $\eta : Id_{DGrph_0} \to UF$ gives the comultiplication $F\eta U: T \to T^2$ and the counit of the adjuction gives $\varepsilon : FU \to Id_{Cat}$, that is, $\varepsilon : T \to Id$.  Now for $\mathbb{A}$ a small category, set $S(\mathbb{A})_n = T^{n+1}(\mathbb{A})$ with face maps $d_i : T^{n+1}(\mathbb{A}) \to T^n(\mathbb{A})$ given by $d_i = T^ {n-i}\varepsilon T^i$, and similarly for the degeneracies which use the comultiplication in an analogous formula. (Note that there are two conventions possible here.  The other will use $d_i = T^i\varepsilon T^{n-i}$. The only effect of such a change is to reverse the direction of certain `arrows' in diagrams later on.  The two simplicial structures are `dual' to each other.)

This $S(\mathbb{A})$ is a simplicial object in $Cat$, $S(\mathbb{A}) : \mathbf{\Delta}^{op} \to Cat$, so does not immediately gives us a simplicially enriched category, however its simplicial set of objects is constant because $U$ and $F$ took note of the identity loops. 

In more detail, let $ob : Cat \to Sets$ be the functor that picks out the set of objects of a small category, then $ob(S(\mathbb{A})) :  \mathbf{\Delta}^{op} \to Sets$ is a constant functor with value the set $ob(\mathbb{A})$ of objects of $\mathbb{A}$. More exactly it is a discrete simplicial set, since all its face and degeneracy maps are bijections.  Using those bijections to identify the possible different sets of objects, yields a constant simplicial set where all the face and degeneracy maps are identity maps, i.e. we do have a \emph{constant} simplicial set.
 
\begin{lemma}
~\\
Let $\mathcal{B} : \mathbf{\Delta}^{op} \to Cat$ be a simplicial object in $Cat$ such that $ob(\mathcal{B})$ is a constant simplicial set with value $B_0$, say.  For each pair $(x,y)\in B_0$, let $$\mathcal{B}(x,y)_{n} = \{\sigma \in \mathcal{B}_{n} | ~  {\rm dom}(\sigma) = x, {\rm codom}(\sigma) = y\},$$
where, of course, ${\rm dom}$ refers to the domain function in $\mathcal{B}_n$ since otherwise ${\rm dom}(\sigma)$ would have no meaning, similarly for $\rm codom$.

(i) The collection $\{\mathcal{B}(x,y)_n |~ n \in \mathbb{N}\}$ has the structure of a simplicial set $\mathcal{B}(x,y)$ with face and degeneracies induced from those of $\mathcal{B}$.

(ii) The composition in each level of $\mathcal{B}$ induces
$$\mathcal{B}(x,y) \times \mathcal{B}(y,z)\to \mathcal{B}(x,z).$$
Similarly the identity map in $\mathcal{B}(x,x)$ is defined as  $id_x$, the identity at $x$ in the category $\mathcal{B}_0$.

(iii) The resulting structure is an $\mathcal{S}$-enriched category, that will also be denoted $\mathcal{B}$.\hfill $\blacksquare$
\end{lemma}

The proof is easy. In particular, this shows that $S(\mathbb{A})$ is a simplicially enriched category.  The description of the simplices in each dimension of $S(\mathbb{A})$ that start at $a$ and end at $b$ is intuitively quite simple.  The arrows in the category, $T(\mathbb{A})$ correspond to strings of symbols representing non-identity arrows in $\mathbb{A}$ itself, those strings being `composable' in as much as the domain of the $i^{th}$ arrow must be the codomain of the $(i-1)^{th}$ one and so on. Because of this we   have:\\
$S(\mathbb{A})_0$ consists exactly of such composable chains of maps in $\mathbb{A}$, none of which is the identity;\\
 $S(\mathbb{A})_1$ consists of such composable chains of maps in $\mathbb{A}$, none of which is the identity, together with a choice of bracketting;\\
 $S(\mathbb{A})_2$ consists of such composable chains of maps in $\mathbb{A}$, none of which is the identity, together with a choice of  two levels of bracketting;\\
and so on.  Face and degeneracy maps remove or insert brackets, but care must be taken when removing innermost brackets as the compositions that can then take place can result in chains with identities which then need removing, see \cite{cordier82}, that is why the comandic description is so much simpler, as it manages all that itself.
\medskip

To understand  $S(\mathbb{A})$ in general it pays to examine the simplest few cases. The key cases are when $\mathbb{A} = [n]$, the ordinal $\{0< \ldots <n\}$ considered as a category in the usual way. The cases $n=0$ and $n=1$ give no surprises. $S[0]$ has one object  0 and $S[0](0,0$ is isomorphic to $\Delta[0]$, as the only simplices are degenerate copies of the identity. $S[1]$ likewise has a trivial simplicial structure, being just the category $[1]$ considered as an $\mathcal{S}$-category. Things do get more interesting at $n = 2$. The key here is the identification of $S[2](0,2)$. There are two non-degenerate strings or paths that lead from 0 to 2, so $S[2](0,2)$ will have two vertices. The bracketted string $((01)(12))$ on removing inner brackets gives $(02)$ and outer brackets, $(01)(12)$ so represents a 1-simplex
$$\xymatrix@+15pt{(01)(12)\ar[r]^{\quad (01)(12))}&(02)}$$
Other simplicial homs are all $\Delta[0]$ or empty.
It thus is possible to visualise $S[2]$ as a copy of $[2]$ with a 2-cell going towards the bottom:
$$
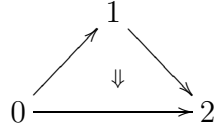
\xymatrix{ &1\ar[dr]&\\
0\ar[rr]\ar[ur]_{\hspace{.5cm}\Downarrow}&&2}$$
The next case $n = 3$ is even more interesting. $S[3](i,j)$ will be (i) empty if $j<i$, (ii) isomorphic to $\Delta[0] $ if $i = j$ or $i = j-1$, (iii) isomorphic to $\Delta[1]$ by the same reasoning as we just saw for $j = i + 2$ and that leaves $S[3](0,3)$. This is a square, $\Delta[1]^2$, as follows:
$$
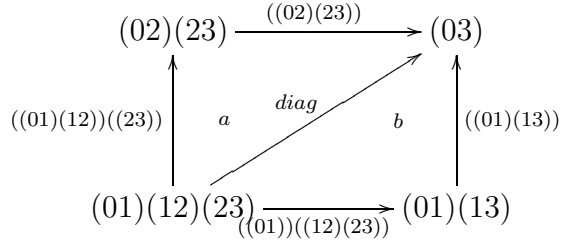
\xymatrix@+25pt{(02)(23)\ar[r]^{((02)(23))}&(03)\\
(01)(12)(23)\ar[u]^{((01)(12))((23))}_{\hspace{.5cm} a}\ar[ur]^{diag}\ar[r]_{((01))((12)(23))}&(01)(13)\ar[u]_{((01)(13))}^{b\hspace{6mm}}
}$$
where the diagonal $diag = ((01)(12)(23))$, $a = (((01)(12))((23)))$ and $b = (((01))((12)(23)))$. (It is instructive to check that this is correct, firstly because I may have slipped up (!) as well as seeing the mechanism in action. Removing the innermost brackets is $d_0$, and so on.)

The case of $S[4]$ is worth doing. I will not draw the diagrams here although aspects of it have implications later, but suggest this as an exercise. As might be expected $S[4](0,4)$ is a cube.

\medskip

\textbf{Remark}

The history of this construction is interesting. A variant of it, but with topologically enriched categories as the end result, is in the work of Boardman and Vogt, \cite{boardmanvogt} and also in Vogt's paper, \cite{vogt73}.  Segal's student Leitch used a similar construction to describe a homotopy commutative
 cube (actually a \emph{homotopy coherent cube}), cf. \cite{leitch}, and this was used by Segal in his famous paper, \cite{segal}, under the name of the `explosion' of $\mathbb{A}$.  All this was still in the topological framework and the link with the comonad resolution was still not in evidence.  Although it seems likely that Kan knew of this link between homotopy coherence and the comonadic resolutions by at least 1980, (cf. \cite{D&K80a}), the construction  does not seem to appear in his work with Dwyer as being linked with coherence until much later. Cordier made the link explicit in \cite{cordier82} and showed how Leitch and Segal's work fitted in to the pattern.  His motivation was for the description of homotopy coherent diagrams of topological spaces.  
Other variants were also apparent in the early work of May on operads, and linked in with Stasheff's work on higher associativity and commutativity `up to homotopy'.

Cordier and Porter, \cite{cordierporter86}, used an analysis of a locally Kan simplicially enriched category involving this construction to prove a generalisation of Vogt's theorem on categories of homotopy coherent diagrams of a given type. (We will return to this aspect a bit later in these notes, but an elementary introduction to this theory can be found in \cite{K&P}.) Finally Bill Dwyer, Dan Kan and Justin Smith, \cite{DKS}, introduced a similar construction for an $\mathbb{A}$ which is an $\mathcal{S}$-category to start with, and motivated it by saying that $\mathcal{S}$-functors with domain this $\mathcal{S}$-category corresponded to \emph{$\infty$-homotopy commutative $\mathbb{A}$-diagrams}, yet they do not seem to be aware of the history of the construction, and do not really justify the claim that it does what they say.  Their viewpoint is however important as, basically, within the setting of Quillen model category structures, this provides a cofibrant replacement construction.  Of course, any other cofibrant replacement could be substituted for it and so would still allow for a description of homotopy coherent diagrams in that context.  This important viewpoint can also be traced to Grothendieck's \emph{Pursuing Stacks}, \cite{stacks}.

The DKS extension of the construction, \cite{DKS}, although simple to do, is often useful and so will be outlined next.  If $\mathbb{A}$ is already a $\mathcal{S}$-category, think of it as a simplicial category, then applying the $\mathcal{S}$-construction to each $\mathbb{A}_n$ will give a bisimplicial category, i.e. a functor $S(\mathbb{A}) : \mathbf{\Delta}^{op} \times \mathbf{\Delta}^{op} \to Cat$. Of this we take the diagonal so the collection of $n$-simplices is $S(\mathbb{A})_{n,n}$, and by noticing that the result has a constant simplicial set of objects, then apply the lemma.

\subsection{The Dwyer-Kan `simplicial groupoid' functor.}
Let $K$ be a simplicial set.  Near the start of simplicial homotopy theory, Kan showed how, if $K$ was reduced (that is, if $K_0$ was a singleton), then the free group functor applied to $K$ in a subtle way, gave a simplicial group whose homotopy groups were those of $K$, with a shift of dimension.  With Dwyer in \cite{D&K84}, he gave the necessary variant of that construction to enable it to apply to the non-reduced case.  This gives a `simplicial groupoid' $G(K)$ as follows: 

The object set of all the groupoids $G(K)_n$ will be in bijective correspondence with the set of vertices $K_0$ of $K$. Explicitly this object set will be written $\{\overline{x} ~|~ x \in K_0\}$.  

The groupoid $G(K)_n$ is generated by edges 
$$\overline{y} : \overline{d_1d_2 \ldots d_{n+1} y} \to \overline{d_0d_2\ldots d_{n+1}y} \quad \mbox{ for } y \in K_{n+1}$$with relations $\overline{s_0x} = id_{\overline{d_1d_2 \ldots d_n x}}$.  Note since these just `kill' some of the generating edges, the resulting groupoid $G(K)_n$ is still a free groupoid. 

Define $\sigma_i\overline{x} = \overline{s_{i+1}x} \quad \mbox{ for } i\geq 0$,
and, for $i > 0$, $\delta_i\overline{x} = \overline{d_{i +1}x}$, but for $i = 0$, $\delta_0\overline{x} = (\overline{d_1x})(\overline{d_0x})^{-1}$.

These definitions yield a simplicial groupoid as is easily checked and, as is clear, its simplicial set of objects  is constant, so it also can be considered as
a  simplicially enriched groupoid, $G(K)$.

(NB. Beware there are several `typos' in the original paper relating to these formulae for the construction and in some of the related material.)

As before it is instructive to compute some examples and we will look at $G(\Delta[2])$ and $G(\Delta[3]$. ) simplicially enriched groupoid are free groupoids in each simplicial dimension, their structure can be clearly seen from the generating graphs.  For instance, $G(\Delta[2])_0$ is the free groupoid on the graph
$$
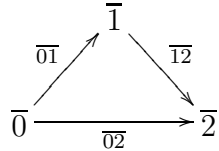
\xymatrix{&\overline{1}\ar[dr]^{\overline{12}}&\\
\overline{0}\ar[rr]_{\overline{02}}\ar[ur]^{\overline{01}}&&\overline{2}}$$		
whilst $G(\Delta[2])_1$ is the free groupoid on the graph
$$
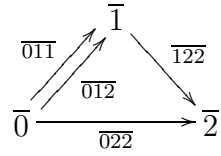
\xymatrix{&\overline{1}\ar[dr]^{\overline{122}}&\\
\overline{0}\ar[rr]_{\overline{022}}\ar[ur]<.5ex>^{\overline{011}}\ar[ur]<-.5ex>_{\overline{012}}&&\overline{2}}$$
Here it is worth noting that $\delta_0(\overline{012}) = (\overline{02}).(\overline{12})^{ -1}$. Higher dimensions do not have any non-degenrate generators.

Again with $G(\Delta[3])$, in dimension 0 we have the free groupoid on the directed graph give by the 1-skelton of $\Delta[3]$. In dimension 2, the generating directed graph is
$$
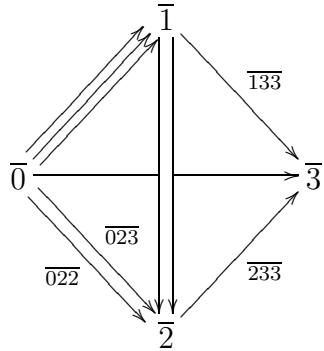
\xymatrix@+20pt{&\overline{1}\ar[dr]^{\overline{133}}\ar[dd]<.5ex>\ar[dd]<-.5ex>&\\
\overline{0}\ar[ur]<.7ex>\ar[ur]\ar[ur]<-.7ex>\ar'[r][rr]\ar[dr]<.5ex>^{\overline{023}}\ar[dr]<-.5ex>_{\overline{022}}&&\overline{3}\\
&\overline{2} \ar[ur]_{\overline{233}}  &   }$$
Here only a few of the arrow labels have been given. Others are easy to provide (but moderately horrible to typeset in a sensible way!). Those from $\overline{0}$ to $\overline{1}$ are $\overline{012}$, $\overline{011}$ and $\overline{013}$; those from $\overline{1}$ to $\overline{2}$ are $\overline{122}$ and $\overline{123}$, and finally from $\overline{0}$ to $\overline{3}$, we have $\overline{033}$.

The next dimension is only a little more complicated. It has extra degenerate arrows such as $\overline{0112}$
and $\overline{0122}$ from $\overline{0}$ to $\overline{1}$ but also between these two vertices has $\overline{0123}$, coming from the non-degenerate three simplex of $\Delta[3]$. The full diagram is easy to draw (and again a bit tricky to typeset in a neat way), and is therefore left `as an exercise'.

\medskip
\textbf{Remarks}

(i) The functor $G$ has a right adjoint $\overline{W}$ and the unit $K\to\overline{W}G(K)$ is a weak equivalence of simplicial sets.  This is part of the result that shows that simplicially enriched groupoids model all homotopy types, for which see the original paper, \cite{D&K84} or the book by Goerss and Jardine, \cite{GoerssJardine}.

(ii) It is tantalising that the definition of $S(\mathbb{A})$ for a category $\mathbb{A}$ and of $G(K)$  for a simplicial set $K$ are very similar yet very different. Why does the twist occur in the $d_0$ of $G(K)$? Why do the source and target maps at each level in $G(K)$ end up at the zeroth and first vertex rather than the seemingly more natural zeroth and $n^{th}$ that occur in $S(\mathbb{A})$?  In fact is there a variant of the $G(K)$ construction that is nearer to the $S(\mathbb{A})$ construction without being merely artificially so?

\section{Structure}
\subsection{The `homotopy' category}
If $\mathcal{C}$ is an $\mathcal{S}$-category, we can form a category
$\pi_0\mathcal{C}$ with the same objects and having $$(\pi_0\mathcal{C})(X,Y) =
\pi_0({\mathcal{C}}(X,Y)).$$  For instance, if $\mathcal{C} =
\mathcal{CW}$, the category of CW-complexes, then $\pi_0\mathcal{CW} = {Ho\mathcal(CW)}$, the corresponding homotopy category.  Similarly we could obtain a groupoid enriched category
  using the fundamental groupoid (cf. Gabriel and Zisman,
  \cite{gabrielzisman}).

One can `do' some elementary homotopy theory in any  $\mathcal{S}$-category,
$\mathcal{C}$, by saying that two maps $f_0, f_1 : X \to Y$ in $\mathcal{C}$ are 
homotopic if there is an $H \in {\mathcal{C}}(X,Y)_1$ with $d_0H =
f_1$, $d_1H = f_0$.

This theory will not be very rich unless at least some low dimensional Kan
conditions are satisfied.  The  $\mathcal{S}$-category, $\mathcal{C}$, is called 
\emph{locally Kan} if each $\mathcal{C}(X,Y)$ is a Kan complex,
\emph{locally weakly Kan} if \ldots , etc. (If you have not met `weak Kan complexes' before, you will soon meet them in earnest!)

The theory is `geometrically' nicer to work with if $\mathcal{C}$ is
\emph{tensored} or \emph{cotensored}.
\subsection{Tensoring and Cotensoring}
\textbf{Tensored}

If for all $K \in \mathcal{S}$, $X, Y,  \in \mathcal{C}$, there is an object
$K\bar{\otimes} X$ in $\mathcal{C}$ such that 
$$\mathcal{C}(K\bar{\otimes}X,Y) \cong
\mathcal{S}(K,\mathcal{C}(X,Y)$$naturally in $K$, $X$ and $Y$, then
$\mathcal{C}$ is said to be \emph{tensored} over $\mathcal{S}$.
  
\medskip

\textbf{Cotensored}

Dually, if we require objects $\bar{\mathcal{C}}(K,Y)$ such that $$\mathcal{C}(X,\bar{\mathcal{C}}(K,Y)) \cong
\mathcal{S}(K,\mathcal{C}(X,Y)$$then we say $\mathcal{C}$ is
\emph{cotensored} over  $\mathcal{S}$.
 
\medskip

To gain some intuitive feeling for these, think of $K\bar{\otimes} X$ as being $K$ product with $X$, and $\bar{\mathcal{C}}(K,Y)$ as the object of functions from $K$ to $Y$.  These words do not, as such, make sense in general, but do tell one the sort of tasks these constructions will be set to do.  They will not be much used explicitly here however.

\begin{proposition}(cf. Kamps and Porter, \cite{K&P})\\
If $\mathcal{C}$ is a locally Kan $\mathcal{S}$-category tensored over
$\mathcal{S}$ then  taking $I\times X = \Delta[1]\bar{\otimes}X$, we get a good
cylinder functor such that for the cofibrations relative to $I$ and weak
equivalences taken to be homotopy equivalences, the category $\mathcal{C}$ has
a cofibration category structure.\hfill $\blacksquare$
\end{proposition}
A cofibration category structure is just one of many variants of the abstract homotopy theory structure introduced to be able to push through homotopy type arguments in particular settings.
There are variants of this result, due to Kamps, see references in \cite{K&P}, where $\mathcal{C}$ is both
tensored and cotensored over $\mathcal{S}$ and the conclusion is that
$\mathcal{C}$ has a Quillen model category structure.  The examples of locally
Kan $\mathcal{S}$-categories include $\mathcal{T}op$, $\mathcal{K}an$,
$\mathcal{G}rpd$ and $\mathcal{C}rs$, but not $\mathcal{C}at$ or $\mathcal{S}$ itself.
\section{Nerves and Homotopy Coherent Nerves.}
\subsection{Kan and weak Kan conditions}
Before we get going on this section, it will be a good idea to bring to the fore the definitions of \emph{Kan complex} and \emph{weak Kan complex} (or \emph{quasi-category}).

 As usual we set $\Delta[n] = {\mathbf \Delta}( - , [n]) \in \mathcal{S}$, then
 for each $i$, $0 \leq i \leq n$, we can form a subsimplicial set, $\Lambda^i[n]$, of
$\Delta[n] $ by discarding the top dimensional $n$-simplex (given by the
identity map on $[n]$) and its $i^{th}$ face.  We must also discard all the
degeneracies of these simplices.  This informal definition does not give a
`picture' of what we have, so we will list the various cases for $n = 2$.
$$\Lambda^0[2] = \xymatrix{
 &1\ar@{ ..}[dr]^{\leftarrow 0^{th} \mbox{~\scriptsize{face missing}}}&\\
                0\ar[ur]\ar[rr]&&2}$$
$$\Lambda^1[2] = \xymatrix{
 &1\ar[dr]^{\hspace{2.3cm}}&\\
                0\ar[ur]\ar@{ ..}[rr]&&2}$$
$$\Lambda^2[2] = \xymatrix{
 &1\ar[dr]^{\hspace{2.3cm}}&\\
                0\ar@{ ..}[ur]\ar[rr]&&2}$$
A map $p: E \rightarrow B$ is a \emph{Kan fibration} if given any $n$, $i$
as above and any $(n,i)$-horn in $E$, i.e. any map $f_1 : \Lambda^i[n]
\rightarrow E$, and $n$-simplex, $f_0 : \Delta[n] \rightarrow B$, such that 
$$
\xymatrix{\Lambda^i[n]\ar[r]^{f_1}\ar[d]_{inc}& E\ar[d]^p\\
\Delta[n]\ar[r]_{f_0}&B}$$
commutes, then there is an $f : \Delta[n] \rightarrow E$ such that $pf = f_0$
and $f.inc = f_1$, i.e. $f$ lifts $f_0$ and extends $f_1$.

A simplicial set, $K$, is a \emph{Kan complex} if the unique map $K
\rightarrow \Delta[0]$ is a Kan fibration.  This is equivalent to saying that
every horn in $K$ has a filler, i.e. any $f_1 : \Lambda^i[n]\rightarrow Y$ extends to an $f : \Delta[n] \rightarrow Y$.  This condition
looks to be purely of a geometric nature but in fact has an important algebraic
flavour; for instance if $f_1 : \Lambda^1[2]\rightarrow Y$ is a horn, it
consists of a diagram
$$
\xymatrix{ & \ar[dr]^b&\\\ar[ur]^a&&}$$
of `composable' arrows in $K$.  If $f$ is a filler, it looks like $$
\xymatrix{ & \ar[dr]^b&\\ \ar[rr]_c\ar[ur]^a_{\quad f}&&}$$
and one can think of the third face $c$ as a composite of $a$ and $b$.  This
`composite' $c$ is not usually uniquely defined by $a$ and $b$, but is `up to
homotopy'. If we write $c = ab$ as a shorthand then if $g_1 :\Lambda^0[2]
\rightarrow K$ is a horn, we think of $g_1$ as being $$
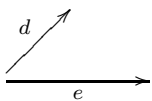
\xymatrix{ & &\\ \ar[ur]^d\ar[rr]_e&&}$$and to find a filler is to find a
diagram  $$
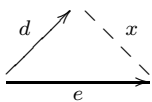
\xymatrix{ & \ar@{--}[dr]^x&\\ \ar[ur]^d\ar[rr]_e&&}$$
and thus to `solve' the equation $dx = e$ for $x$ in terms of $d$ and
$e$. It thus requires in general some approximate inverse for $d$, in fact, taking $e$ to be a degenerate 1-simplex puts one in exactly such a position. 

In many useful cases, we do not always have inverses  and so want to discard any requirement that would imply they always exist.  This leads to the weaker form of the Kan condition in which in each dimension no requirement is made for the existence of fillers on horns that miss out the zeroth or last faces.  More exactly:

\medskip

\textbf{Definition}

 A simplicial set ${\bf K}$ is \emph{a weak Kan complex} or \emph{quasi-category} if for any $n$ and $0< k < n$, any $(n,k)$-horn in $K$ has a filler.

\medskip

\textbf{Remarks}

(i) Joyal, \cite{joyal}, uses the term \emph{inner horn} for any $(n,k)$-horn in $K$ with $0< k < n$. The two remaining cases are then conveniently called \emph{outer horns}.

(ii) For any space $X$, its singular complex, $Sing(X)$ is given by $Sing(X)_n = Top(\Delta^n,X)$ with the well known face and degeneracy maps.  This simplicial set is always a Kan complex as is the related $\mathcal{T}op(X,Y)$ as mentioned above. 

\medskip

\subsection{Nerves}

The categorical analogue of the singular complex is, of course, the nerve:  if $\mathbb{C}$ is a category, its \emph{nerve}, $Ner(\mathbb{C})$, is the simplicial set with $Ner(\mathbb{C})_n = Cat([n],\mathbb{C})$, where $[n]$ is the category associated to the finite ordinal $[n] = \{ 0 < 1< \ldots < n\}$.  The face and degeneracy maps are the obvious ones using the composition and identities in $\mathbb{C}$. 

The following is well known and easy to prove.

\begin{lemma}
~\\
(i) $Ner(\mathbb{C})$ is always weakly Kan. \\
(ii) $Ner(\mathbb{C})$ is Kan if and only if $\mathbb{C}$ is a groupoid.\hfill $\blacksquare$
\end{lemma}

Of course more is true. Not only does any inner horn in $Ner(\mathbb{C})$ have a filler, it has exactly one filler.  To express this with maximum force the following idea, often attributed to Graeme Segal or to Grothendieck, is very useful. 

Let $p>0$, and consider the increasing maps $e_i : [1] \to [p]$ given by $e_i(0) = i$ and $e_i(1) = i+1$.  For any simplicial set $A$ considered as a functor $A : \mathbf{\Delta}^{op} \to Sets$, we can evaluate $A$ on these $e_i$ and, noting that $e_i(1) = e_{i+1}(0)$, we get a family of functions $A_p \to A_1$, which yield a cone diagram, for instance, for $p =3$:
$$
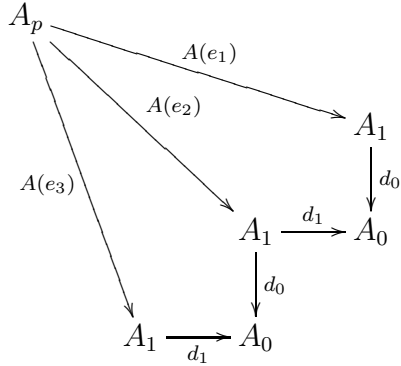
\xymatrix{A_p \ar[drrr]^{A(e_1)}\ar[ddrr]^{A(e_2)}\ar[dddr]_{A(e_3)}&&&\\
&&&A_1\ar[d]^{d_0}\\
&&A_1\ar[r]^{d_1}\ar[d]^{d_0}&A_0\\
&A_1\ar[r]_{d_1}&A_0&&}$$
and in general, thus yield a map
$$\delta[p]: A_p \to A_1\times_{A_0}A_1  \times_{A_0}\ldots \times_{A_0}A_1.$$
The maps, $\delta[p]$, have been called the \emph{Segal maps} and will recur throughout the rest of these notes.

\begin{lemma}
~\\
If  $A = Ner(\mathbb{C})$ for some small category $\mathbb{C}$, then for $A$, the Segal maps are bijections.
\end{lemma}
\textbf{Proof}

A simplex $\sigma \in Ner(\mathbb{C})_p$ corresponds uniquely to a composable $p$-chain of arrows in $\mathbb{C}$, and hence exactly to its image under the relevant Segal map. \hfill$\blacksquare$

\medskip

Better than this is true:

\begin{proposition}\label{GrotSegal}
~\\
If $A$ is a simplicial set such that the Segal maps are bijections then there is a category structure on the directed graph $$\xymatrix{A_1 \ar[r]<1ex> \ar[r]&A_0\ar[l]<1ex>}.$$making it a category whose nerve is isomorphic to the given $A$.
\end{proposition}
\textbf{Proof}

 To get composition you use $$A_1\times_{A_0}A_1 \stackrel{\cong}{\to} A_2\stackrel{d_1}{\to}A_1.$$
Associativity is given by $A_3$.  The other laws are easy, and illuminating, to check. \hfill$\blacksquare$
\medskip

The condition `Segal maps are a bijection' is closely related to notions of `thinness' as used by Brown and Higgins in the study of crossed complexes and their relationship to $\omega$-groupoids, see, for instance, \cite{rb&pjh1987}, and also to Duskin's `hypergroupoid' condition, \cite{duskin}.

Another result that is sometimes useful, is a refinement of `groupoids give Kan complexes':.  The proof is `the same':
\begin{lemma}
~\\
Let $A = Ner(\mathbb{C})$, the nerve of a category $\mathbb{C}$.

(i) Any $(n,0)$-horn
$$f : \Lambda^0[n]\rightarrow A$$
for which $f(01)$ is an isomorphism has a filler. Similarly any $(n,n)$-horn $g : \Lambda^n[n]\rightarrow A$ for which $g(n-1~n)$ is an isomorphism, has a filler.

(ii) Suppose $f$ is a morphism in $\mathbb{C}$ with the property that for any $n$, any $(n,0)$-horn $\phi : \Lambda^0[n]\rightarrow A$ having $f$ in the $(0,1)$ position, has a filler, then $f$ is an isomorphism.  (Similarly with $(n,0)$ replaced by $(n,n)$ with the obvious changes.)
\hfill$\blacksquare$
\end{lemma}

\medskip

\textbf{Remark}

Joyal in \cite{joyal} suggested that the name `weak Kan complex', as introduced by  Boardman and Vogt, \cite{boardmanvogt}, could be changed to that of `quasi-category' to stress the analogy with categories \emph{per se} as `\emph{Most concepts and results of category theory can be extended to quasi-categories}', \cite{joyal}.

It would have been nice to have explored Joyal's work on quasi-categories more fully, e.g. \cite{joyal}, but time did not allow it.  The following few sections just skate the surface of the theory.

\subsection{Quasi-categories}
Categories yield quasi-categories via the nerve construction.  Quasi-categories yield categories by a `fundamental category' construction that is left adjoint to nerve. This can be constructed using the free category generated by the 1-skeleton of $A$, and then factoring out by a congruence generated by the basic relations : $gf \equiv h$, one for each commuting 1-sphere $(g,h,f)$ in $A$.  By a \emph{1-sphere} is meant a map $a : \partial \Delta[2] \to A$, thus giving three faces, $(a_0,a_1,a_2)$ linked in the obvious way. The 1-sphere is said to be \emph{commuting} if there is a 2-simplex, $b\in A_2$, such that $a_i = d_ib$ for $i = 0,1,2$.     

This `fundamental category' functor also has a very neat description due to Boardman and Vogt.  (The treatment here is adapted from \cite{joyal}.)

We assume given a quasi-category $A$.  Write $gf \sim h$ if $(g,h,f)$ is a commuting 1-sphere.  Let $x,y \in A_0$ and let $A_1(x,y) = \{f\in A_1 ~|~ x = d_1f, y = d_0f\}$.  If $f,g \in A_1(x,y)$, then, suggestively writing $s_0x = 1_x$, 
\begin{lemma}~\\
The four relations $f1_x \sim g$, $g1_x\sim f$, $1_yf\sim g$ and $1_yg\sim f$ are equivalent.\hfill$\blacksquare$
\end{lemma}

The proof is easy.

We will say $f\simeq g$ if any of these is satisfied and call $\simeq$, the \emph{homotopy relation}. It is an equivalence relation on $A_1(x,y)$.  Set $ho\,A_1(x,y) = A_1(x,y)/\simeq$.

If $f \in A_1(x,y)$, $g \in A_1(y,z)$ and $h\in A_1(x,z)$, then the relation $gf \sim h$ induces a map:
$$ho\,A_1(x,y)\times ho\,A_1(y,z) \to ho\,A_1(x,z).
$$
\begin{proposition}~\\
The maps $$ho\,A_1(x,y)\times ho\,A_1(y,z) \to ho\,A_1(x,z)$$give a composition law for a category, $ho\,A$, the homotopy category of $A$.\hfill$\blacksquare$
 \end{proposition}

Of course, $ho\,A$ is the fundamental category of $A$ up to natural isomorphism. From previous comments we have:

\begin{corollary}~\\
A quasi-category $A$ is a Kan complex if and only if $ho\,A$ is a groupoid.\hfill$\blacksquare$
\end{corollary}
\medskip

\subsection{Homotopy coherent nerves}

Before introducing this topic, recall some of the intuition behind homotopy coherent (h.c.) diagrams. (Again there is an overview of this theory in \cite{cubo} and a thorough introduction in \cite{K&P}.)

\medskip

\textbf{Examples of h.c. diagrams in $Top$.}

1)  A diagram indexed by the small category, $[2]$,
$$
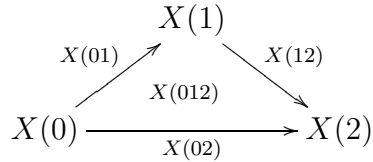
\xymatrix{ &X(1)\ar[dr]^{X(12)}&\\
X(0)\ar[rr]_{X(02)}\ar[ur]^{X(01)}_{\hspace*{.5cm} X(012)}&&X(2)}$$
is h.c. if there is specified a homotopy 
$$X(012) : I\times X(0) \to X(2),$$
$$X(012) : X(02) \simeq X(12)X(01).$$

\medskip

2) For a diagram indexed by $[3]$: Draw a 3-simplex, marking the vertices
$X(0)$, \ldots, $X(3)$, the edges $X(ij)$, etc., the faces $X(ijk)$, etc.  The 
homotopies $X(ijk)$ fit together to make the sides of a square
$$
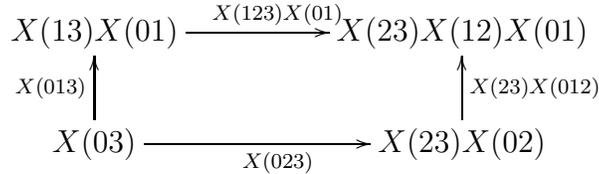
\xymatrix{X(13)X(01)\ar[rr]^{X(123)X(01)} && X(23)X(12)X(01)\\
X(03)\ar[u]^{X(013)}\ar[rr]_{X(023)}&&X(23)X(02)\ar[u]_{X(23)X(012)}}$$
and the
diagram is made h.c. by specifying a second level homotopy 
$$X(0123) : I^2\times X(0)\to X(3)$$filling this square.

\medskip

These can be continued for larger $[n]$.  Of course, this is not how the theory is formally specified, but it
provides some understanding of the basic idea. 

 The theory was initially developed 
by Vogt, \cite{vogt73}, following methods introduced with Boardman,
\cite{boardmanvogt} (see also the references in that source for other
earlier papers on  the area). Cordier \cite{cordier82} provides a simple
$\mathcal{S}$-category theory way of working with h.c. diagrams and hence
released an `arsenal' of categorical tools for working with h.c. diagrams.  Some of that is worked out in the papers, \cite{C&P88,C&P90,C&P96,C&P97}

\medskip

\textbf{Some Results}

(i) If $X : {\mathbb{A}}\to {\mathcal{T}op}$ is a commutative diagram and we
replace some of the $X(a)$ by homotopy equivalent $Y(a)$ with specified
homotopy equivalence data:
$$f(a) : X(a) \to Y(a), \quad g(a) : Y(a) \to X(a)$$
$$H(a) : g(a)f(a) \simeq Id, \quad K(a) : f(a)g(a) \simeq Id,$$
then we can combine these data into the construction of a h.c. diagram $Y$
based on the objects $Y(a)$ and homotopy coherent maps $$f : X\to Y, \quad g : 
Y \to X, \mbox { etc.,}$$making $X$ and $Y$ homotopy equivalent as h.c. diagrams.

(This is `really' a result about quasi-categories, see \cite{joyal}.)
\medskip

(ii) {Vogt}, \cite{vogt73}.

If ${\mathbb{A}}$ is a small category, there is a category ${Coh(\mathbb{A},\mathcal{T}op)}$
  of h.c. diagrams and homotopy classes of h.c. maps between them.  Moreover
  there is an equivalence of categories
$${Coh(\mathbb{A},\mathcal{T}op)} \stackrel{\simeq}{\to} {Ho(\mathcal{T}op^\mathbb{A})}$$
This was extended replacing $\mathcal{T}op$ by  a general locally Kan simplicially enriched complete category, $\mathcal{B}$,  in \cite{cordierporter86}.
\medskip

(iii) {Cordier (1980)}, \cite{cordier82}.

Given ${\mathbb{A}}$, a small category, then the $\mathcal{S}$-category
${S(\mathbb{A})}$ is such that a
h.c. diagram of type ${\mathbb{A}}$ in ${\mathcal{T}op}$ is given precisely by an 
$\mathcal{S}$-functor
$$F : {S(\mathbb{A})} \to {\mathcal{T}op}$$
This suggested the extension of h.c. diagrams to other contexts such as a
general locally Kan $\mathcal{S}$-category, $\mathcal{B}$ and suggests the definition of homotopy coherent diagram in a $\mathcal{S}$-category and thus a h.c. nerve of an $\mathcal{S}$-category.

\vspace{.5cm}
 
\textbf{Definition}  ({Cordier (1980)}, \cite{cordier82}, based on earlier ideas of Vogt, and Boardman-Vogt.)

Given a simplicially enriched category $\mathcal{B}$, the \emph{homotopy coherent nerve} of $\mathcal{B}$, denoted $Ner_{h.c.}(\mathcal{B})$, is the simplicial `set' with 
$$Ner_{h.c.}(\mathcal{B})_n = \mathcal{S}\!-\!Cat(S[n],\mathcal{B}).$$

\medskip

To understand simple h.c. diagrams and thus $Ner_{h.c.}(\mathcal{B})$, we will unpack the definition of homotopy coherence.  

The first thing to note is that for any $n$ and $0\leq i< j\leq n$, $S[n](i,j) \cong \Delta[1]^{j-i-1}$, the $(j-i-1)$-cube given by the product of $j-i-1 $ copies of $\Delta[1]$.  Thus we can reduce the higher homotopy data to being just that, maps from higher dimensional cubes.

Next some notation:

Given simplicial maps
$$f_1: K_1 \to \mathcal{B}(x,y),$$
$$f_2: K_2 \to \mathcal{B}(y,z),$$
we will denote the composite
$$K_1 \times K_2 \to \mathcal{B}(x,y)\times \mathcal{B}(y,z) \stackrel{c}{\to} \mathcal{B}(x,z)$$
just by $f_2.f_1$ or $f_2f_1$.  (We already have seen this in the h.c. diagram above for $\mathbb{A} = [3]$.  $X(123)X(01)$ is actually $X(123)(I \times X(01) )$, whilst $X(23)X(012)$ is exactly what it states.)

Suppose now that we have the h.c. diagram $F : S(\mathbb{A}) \to \mathcal{B}$.  This is an $\mathcal{S}$-functor and so:\\
to each object $a$ of $\mathbb{A}$, it assigns an object $F(a)$ of $\mathcal{B}$;\\
for each string of composable maps in $\mathbb{A}$,
$$\sigma = (f_0, \ldots, f_n)$$
starting at $a$ and ending at $b$, a simplicial map
$$F(\sigma) : S(\mathbb{A})(0,n+1) \to \mathcal{B}(F(a), F(b)),$$
that is, a higher homotopy
$$F(\sigma) : \Delta[1]^n \to \mathcal{B}(F(a), F(b)),$$
such that

(i) if $f_0 = id$, $F(\sigma) = F(\partial_0\sigma)(proj \times \Delta[1]^{n-1})$

(ii) if $f_i = id$, $0< i < n$
$$F(\sigma) = F(\partial_i\sigma(.(I^i \times m \times I^{n-i}),$$ where $m : I^2 \to I$ is the multiplicative structure on $I = \Delta[1]$ by the `max' function on $\{0,1\}$;

(iii) if $f_n = id$, $F(\sigma) = F(\partial_n \sigma)$;

(iv)$_{i}$  $F(\sigma)|(I^{i-1}\times \{0\} \times I^{n-i}) = F(\partial_i\sigma), 1 \leq i \leq n-1$;

(v)$_{i}$  $F(\sigma)|( I^{i-1}\times \{1\} \times I^{n-i}) = F(\sigma^\prime_i) . F(\sigma_i)$, where $\sigma_i = (f_0, \ldots, f_{i-1}$ 
and $\sigma^\prime = (f_i, \ldots, f_n)$.  We have used $\partial_i$ for the face operators in the nerve of $\mathbb{A}$.

\medskip

The specification of such a homotopy coherent diagram can be split into two parts: \\
(a) specification of certain homotopy coherent \emph{simplices}, i.e. elements in $Ner_{h.c.}(\mathcal{B})$; \\
and\\
(b) specification, via a simplicial mapping from $Ner(\mathbb{A})$ to $Ner_{h.c.}(\mathcal{B})$, of how these individual parts (from (a)) of the diagram are glued together.

The second part of this is easy   as it amounts to a simplicial map $Ner(\mathbb{A}) \to Ner_{h.c.}(\mathcal{B})$, and so we are left with the first part.  The following theorem was proved by Cordier and myself, but the idea was essentially in Boardman and Vogt's lecture notes, like so much else!

\begin{theorem} (\cite{cordierporter86})\\
If $\mathcal{B}$ is a locally Kan $\mathcal{S}$-category then $Ner_{h.c.}(\mathcal{B})$ is a quasi-category.\hfill$\blacksquare$
\end{theorem}
\medskip

It seems to be the case that if $\mathcal{B}$ is only locally weakly Kan, then $Ner_{h.c.}(\mathcal{B})$ need not be a quasi-category.

The proof of the theorem is in the paper, \cite{cordierporter86} and is not too complex. The essential feature is that the very definition (unpacked version) of homotopy coherent diagram makes it clear that parts of the data have to be composed together, (recall the composition of simplicial maps
$$f_1: K_1 \to \mathcal{B}(x,y),$$
$$f_2: K_2 \to \mathcal{B}(y,z),$$
above and how important that was in the unpacked definition).

\medskip

We thus have that a homotopy coherent diagram `is' a simplicial map $F : Ner(\mathbb{A}) \to Ner_{h.c.}(\mathcal{B})$  and that $Ner_{h.c.}(\mathcal{B})$ is a quasi-category.  Of course, the usual proof that, if $X$ and $Y$ are simplicial sets, and $Y$ is Kan, then $\mathcal{S}(X,Y)$ is Kan as well, extends to having $Y$ a quasi-category and the result being a quasi-category.  Earlier we referred to $Coh(\mathbb{A},\mathcal{B})$  in connection with Vogt's theorem.
The neat way of introducing this is as $ho\, \mathcal{S}(Ner(\mathbb{A}), Ner_{h.c.}(\mathcal{B}))$, the fundamental category of the function quasi-category.  In fact, this is essentially the way Vogt first described it.

\vspace{.5cm}

Before we leave homotopy coherence, there is a point that is worth noting for the links with algebraic and categorical models for homotopy types.  The $\mathcal{S}$-categories, $S[n]$, contain a lot of the information needed for the construction of such models.  A good example of this is the interchange law and its links with Gray categories and Gray groupoids.  

Consider $S[4]$.  The important information is in the simplicial set $S[4](0,4)$.  This is a 3-cube, so is still reasonably easy to visualise.  Here it is. The notation is not intended to be completely consistent with earlier uses but is meant to be more self explanatory.
$$
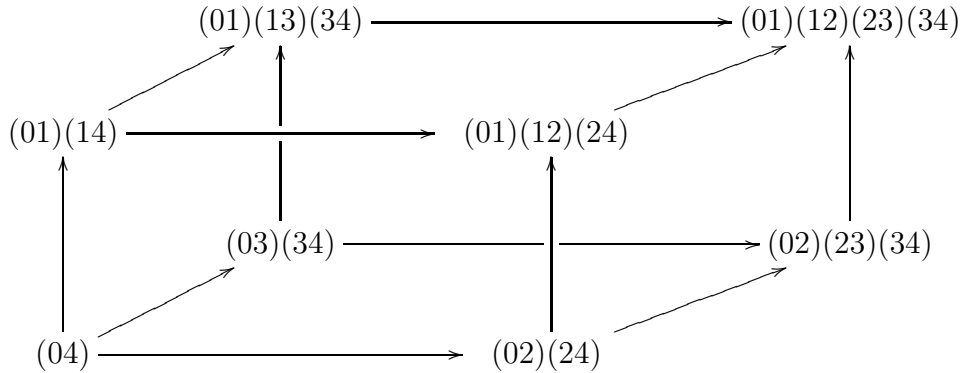
\xymatrix{& (01)(13)(34)\ar[rr]&&(01)(12)(23)(34)\\
(01)(14)\ar[rr]\ar[ur]&&~~(01)(12)(24)\quad\ar[ur]&\\
& (03)(34)\ar'[r][rr]\ar'[u][uu]&&(02)(23)(34\ar[uu])\\
(04)\ar[rr]\ar[ur]\ar[uu]&&~~(02)(24)\quad\ar[ur]\ar[uu]&}$$
This looks mysterious! A 4-simplex has 5 vertices, and hence 5 tetrahedral faces.  Each of the 5 tetrahedral faces will contribute a square to the above diagram, yet a cube has 6 square faces!  (Things get `worse' in $S[5](0,5)$, which is a 4-cube, so has  8 cubes as its faces, but $\Delta[5]$ has only 6 faces.)  Back to the extra face,  this is
$$
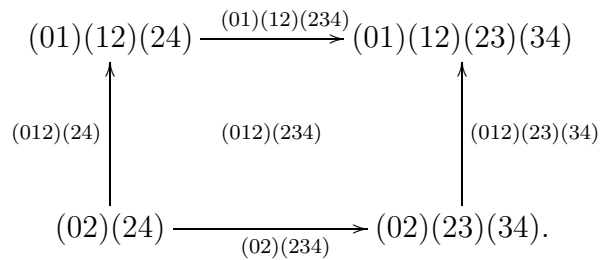
\xymatrix{(01)(12)(24)\ar[rr]^{(01)(12)(234)}&&(01)(12)(23)(34)\\
&&\\
(02)(24)\ar[rr]_{(02)(234)}\ar[uu]^{(012)(24)}&\ar@{}[uu]|{(012)(234)}&(02)(23)(34)\ar[uu]_{(012)(23)(34)}.}$$
The arrow $(012) : (02) \rightarrow (01)(12)$ will, in a homotopy coherent diagram, make its appearence as the  homotopy,$$X(012) : I\times X(0) \to X(2),$$
$$X(012) : X(02) \simeq X(12)X(01),$$
thus this square implies that the homotopies $X(012)$ and $X(234)$ interact minimally. Drawing them as 2-cells in the usual way, the square we have above is the interchange square and the interchange law will hold in this system provided this square is, in some sense, commutative.  In models for homotopy $n$-types for $n \geq 3$, these interchange squares give part of the pairing structure between different levels of the model.  They are there in, say, the Conduch\'e model (2-crossed modules) as the Peiffer lifting, (cf. Conduch\'e, \cite{conduche}) and in the Loday model, (crossed squares, cf. \cite{loday}), as the $h$-map.  In a general dimension, $n$, there will be pairings like this for any splitting of $\{0,1, \ldots ,n\}$ of the form $\{0.1. \ldots, k\}$ and $\{k, \ldots, n\}$.

\section{Dwyer-Kan Hammock Localisation: more simplicially enriched categories.}

In his original contribution \cite{quillen} to abstract homotopy theory, Quillen introduced the notion of a \emph{model category}.  Such a context is a category, $C$, together with three classes of maps: weak equivalences, $W = C_{w.e.}$; fibrations, $fib = C_{fib}$; and cofibrations, $cofib = C_{cofib}$, satisfying certain axioms so as to give a general framework for `doing homotopy theory'.  One of the constructions he used was a categorical localisation already well known from Gabriel's thesis and the work of the French school of algebraic geometers, (Grothendieck, Verdier, etc.) and concurrently with the publication of \cite{quillen}, studied in some depth by Gabriel and Zisman, \cite{gabrielzisman}. The main point was that the analogues of homotopy equivalences, in important instances of homotopical or homological algebra, were only `weak equivalences' so whilst with a homotopy equivalence between two spaces, you are given two maps, one in each direction, plus of course some homotopies, when you have, for instance, a quasi-isomorphism between two chain complexes, you only had one map in one direction: $f: C\to D$ together with the knowledge that $f_*: H_*(C) \to H_*(D)$ was an isomorphism.  The partial solution was to go to the `homotopy category' by formally inverting the weak equivalences, thus getting formal maps going in the opposite direction! (This  may look like cheating, but really is no worse than introducing fractions into the integers, so as to be able to solve certain equations, and of course the detailed construction is closely related!)  We thus end up with a category $C[W^{-1}]$.

This construction is very useful, but this homotopy category does \emph{not} capture the higher order homotopy information implicit in $C$. For instance, the problem of the `best' way to handle homotopy limits and colimits, and more generally derived Kan extensions,  in a model category setting is still central to much of the work on abstract homotopy theory, (cf. \emph{Les D\'erivateurs}, by George Maltsiniotis, \cite{malt1} see also \cite{malt2}, Denis-Charles Cisinski's thesis, and subsequent work, (cf. \cite{cisinski1,cisinski2} and related papers), the resum\'e of Thomason's note books published by Chuck Weibel, \cite{weibel} and Carlos Simpson's, \cite{carlos9708010}, for example).
In a series of articles \cite{D&K80a,D&K80b,D&K80c} published in 1980, Dwyer and Kan proposed a neat solution to this problem,  simplicial localisations.  We will limit ourselves to one of the two versions here, the hammock localisation.
\subsection{Hammocks}
Given a category $C$, and a subcategory $W$, having the same class of objects, construct a $\mathcal{S}$-category $L^H(C,W)$ or $L^HC$ for short, the \emph{hammock localisation of $C$ with respect to $W$}, as follows:

The objects of $L^HC$ are the same as those of $C$

Given two objects $X$ and $Y$, the $k$-simplices of $L^HC(X,Y)$ will be the ``reduced hammocks of width $k$ and any length'' between $X$ and $Y$.
Such a thing is a commutative diagram of form
$$
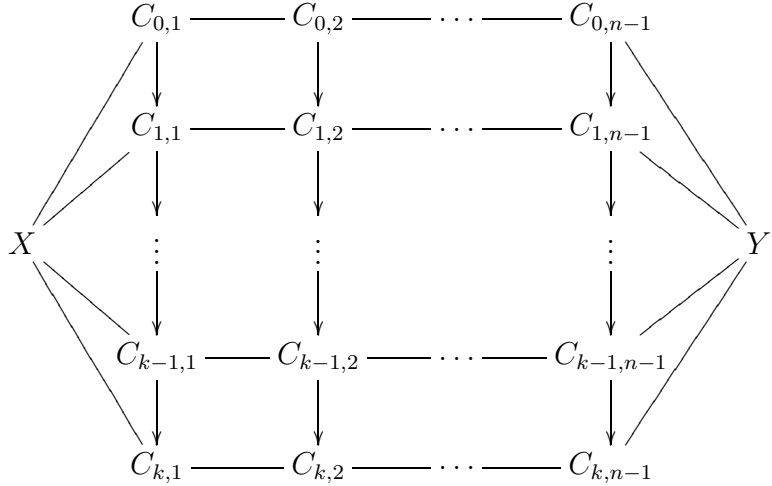
\xymatrix{&C_{0,1}\ar[d]\ar@{-}[r]&C_{0,2}\ar[d]\ar@{-}[r]&\ldots\ar@{-}[r]&C_{0,n-1}\ar[d]\ar@{-}[ddr]&\\
		&C_{1,1}\ar[d]\ar@{-}[r]&C_{1,2}\ar[d]\ar@{-}[r]&\ldots\ar@{-}[r]&C_{1,n-1}\ar[d]\ar@{-}[dr]&\\
	X\ar@{-}[dr]\ar@{-}[ddr]\ar@{-}[uur]\ar@{-}[ur]& \vdots\ar[d] & \vdots\ar[d]&&\vdots\ar[d]& Y\\	
	&C_{k-1,1}\ar[d]\ar@{-}[r]&C_{k-1,2}\ar[d]\ar@{-}[r]&\ldots\ar@{-}[r]&C_{k-1,n-1}\ar[d]\ar@{-}[ur]&\\	
	&C_{k,1}\ar@{-}[r]&C_{k,2}\ar@{-}[r]&\ldots\ar@{-}[r]&C_{k,n-1}\ar@{-}[uur]&}$$
in which\\
(i) the length $n$ of the hammock can be any integer $\geq 0$,\\
(ii)  all the vertical maps are in $W$,\\
(iii) in each column of horizontal maps, all maps go in the same direction; if they go left, then they have to be in $W$;\\
plus two reduction conditions,\\
(iv) the maps in adjacent columns go in different directions,\\
and\\
(v) no column contains only identity maps.

(If in manipulating hammocks, these last two conditions become violated. then it is simple to reduce the hammock by, for example, composing adjacent columns if they point in the same direction or by removing a column of identities.  Repeated use of the reductions may be needed.  One reduction may create a need for another one. It is often useful to work with unreduced hammocks and then to reduce.)

The face and degeneracy maps are defined in the obvious way, (remember the vertices of such a simplex are the `zigzags' from $X$ to $Y$), however they may result in a non-reduced hammock. 

Composition is by concatenation followed by reduction:
$$L^HC(X,Y) \times L^HC(Y,Z)\to L^HC(X,Z),$$
expanding the intervening $Y$ node into a vertical line with identities and then reducing if need be.

Each $L^HC(X,Y)$ is the direct limit of nerves of small categories in an obvious way, i.e. increasing the length $n$ of the hammocks, and so is itself  a quasi-category. 

\subsection{Hammocks in the presence of a calculus of left fractions.}
 If the pair $(C,W)$ satisfies  any of the usual `calculus of fractions' type conditions, then the homotopy type of those nerves already stabilises early on in the process (i.e. for small $n$). The argument given in \cite{D&K80b} is indirect, so let us briefly see why one of these claims is true.
Suppose that $(C,W)$ satisfies a calculus of left fractions, then \\
(i) whenever there is a diagram $X^\prime \stackrel{u}{\leftarrow} X \stackrel{f}{\to} Y$ in $C$ with $u \in W$, then there is a diagram $X^\prime \stackrel{f^\prime}{\to}Y^\prime \stackrel{v}{\leftarrow} Y$ so that $v\in W$ and $vf = f^\prime u$, \\and similarly\\
 (ii) if $f,g: X\to Y \in C$ and $u : X\to X^\prime \in W$ is such that $fu = gu$, then there is a $v\in W$ such that $vf = vg$. 
\\
(By this means any word in arrows of $C$ and $W^{-1}$ can be rewritten to get all the occurrences of arrows from $W^{-1}$ to the left of those `ordinary' arrows from $C$. Each of the two substrings can then be composed to reduce the word to one of the form $w^{-1}c$, i.e. a left fraction.)  To understand how this reacts with hammocks, consider a simple case where the chosen vertex of the hammock $L^HC(X,Y)$ is simply
$$\xymatrix{X&C\ar[l]_w\ar[r]^c& Y &Y\ar[l]_{id}}$$provide 
with $w\in W$.  We construct a new diagram, using the left fractions rule (i), giving a 1-simplex with the given vertex at one end:
$$\xymatrix{X\ar@{=}[d]&C\ar[l]_w\ar[r]^c\ar[d]^w&Y\ar[d]^{w^\prime} &Y\ar[l]_{id}\ar@{=}[d]\\
	        X	               &X\ar[l]_{id}\ar[r]_{c^\prime}    & C^\prime&Y\ar[l]_{w^\prime} },$$
so was homotopic to a `left biased' hammock $(w^\prime)^{-1}c^\prime$.
 
Of course, if the length of the hammock had been greater then the chain of `moves' to link it to the `left biased ' form would be longer.  Again of course, although combinatorially feasible a detailed proof that the left baissed hammocks with   vertices of the form $$X\to C \leftarrow Y$$ provide a deformation retract of $L^HC(X,Y)$ is technically quite messy.  Even with a better knowledge of what the $L^HC(X,Y)$ looks like, there is still the problem of composition.  Two left biased hammocks compose by concatenation to give a more general form of hammock that then gets reduced by the left fractions rules, but these rules do \emph{not} give a normal form for the composite.  Much as in the composite of arrows in a quasi-category, the composite here is only defined up to homotopy.  

Suppose we let $L^1(X,Y)$ be the simplicial set of such left biased hammocks, then it is a deformation retract of $L^HC(X,Y)$. After composition we reduce to get a diagram
$$\xymatrix{L^1(X,Y)\times L^1(Y,Z) \ar[r]\ar@{^{(}->}[d]_{\simeq}\ar[dr]^{concat} &L^1(X,Z)\ar@{^{(}->}[d]^\simeq \\
		L^HC(X,Y)\times L^HC(Y,Z) \ar[r] &L^HC(X,Z)\ar[u]<1ex>^{reduce}}$$ 
This looks as if it should work well, but if we look at the associativity axiom, it is represented by a commutative diagram, and we have replaced each of the nodes of that diagram by a homotopy equivalent object, so we risk getting a homotopy coherent diagram, not a commutative one.
This is happening inside $L^HC$,  so this does not matter so much.  Although attempting to cut down the size of the `hom-sets' does allow us more control over some aspects of the situation, it also has its downside.

The solution is to study the homotopy theory of $\mathcal{S}$-categories as such.  This will lead us back towards the Segal maps as well as continuing to interact with homotopy coherence.

For a short time, for the purpose of exposition, we will restrict ourselves to small $\mathcal{S}$-categories with a fixed set of objects, $O$, say, and $\mathcal{S}$-functors will be the identity on objects. We will denote the category of such things by $S\!-\!Cat/O$. (The material here is adapted from \cite{DKS}.)  This category has a closed simplicial model category structure in which the simplicial structure is more or less obvious, in which a map $D\to D^\prime$ is a weak equivalence (resp.  a fibration), whenever, for every pair of objects, $x,y \in O$, the restricted map
$$D(x,y)\to D^\prime(x,y)$$
is a weak equivalence (resp. fibration).  (Note, that several of the constructions we have been looking at gave us weak equivalences in this sense, for instance, $S(\mathbb{A})\to \mathbb{A}$ is one such and that the fibrant objects are the locally Kan $\mathcal{S}$-categories over $O$).

Now as we know any of the categories $S\!-\!Cat/O$ form subcategories of the category of simplicial categories, $Cat^{\mathbf{\Delta}^{op}}$.  This category also has a closed simplicial model category structure and the nerve and categorical realisation functors induce an equivalence of homotopy categories  (even of the simplicial localisations if you want) between $Cat^{\mathbf{\Delta}^{op}}$ and the category of bisimplicial sets $S^{\mathbf{\Delta}^{op}}$.  Within $Cat^{\mathbf{\Delta}^{op}}$ we are used to considering $S\!-\!Cat$ as a full subcategory. Related to the problem of reducing the size of the $L^HC(X,Y)$s is the question of determining the result of restricting the induced nerve functor to $S\!-\!Cat$. The solution is rather surprising:

Consider the full subcategory of $S^{\mathbf{\Delta}^{op}}$ determined by those objects $X$ such that (i) $X[0]$ is a discrete simplicial set (cf. the condition on the object simplicial set in an $\mathcal{S}$-category); \\
and\\
(ii) for every integer $p \geq 2$, the Segal map $$\delta[p]: X[p] \to X[1]\times_{X[0]}X[1]\times_{X[0]}\ldots \times_{X[0]}X[1]$$
is a weak equivalence of simplical sets.

For reasons that will become clearer later, we will call these objects \emph{Segal categories} or sometimes \emph{Segal 1-categories}. Of course, there is a notion of Segal 0-categories, but these are just nerves of ordinary categories.  We will denote the category of these Segal 1-categories by $Segal\!-\!Cat$.  The result of Dwyer, Kan and Smith, \cite{DKS}, is that the nerve from $Cat^{\mathbf{\Delta}^{op}}$ to $S^{\mathbf{\Delta}^{op}}$, restricts to given an equivalence of homotopy categories between $S\!-\!Cat$ and $Segal\!-\!Cat$.  In particular this says that any Segal category is weakly equivalent to a bisimplicial set that is a nerve of a simplicially enriched category.  Segal categories are weakened simplicial versions of the algebraic structures given by the categorical axioms, so this is in many ways a coherence theorem for Segal categories rather like the coherence theorems for bicategories, etc. 

\section{$\mathbf{\Gamma}$-spaces, $\mathbf{\Gamma}$-categories and Segal categories.}

The reason that Segal categories arise as they do is best sought in the paper \cite{segal} by Segal, although it is not there but rather in \cite{DKS} that they were introduced, but not named as such. (In fact their first naming seems to be in Simpson's \cite{carlos9710011}.) In \cite{segal}, one of the main aims was to get `up-to-homotopy' models for algebraic structures so as to be able to iterate classifying space constructions, to form spectra 
for studying corresponding cohomology theories and to help `delooping' spaces where appropriate.  Various approaches had been tried, notably that of Boardman and Vogt, \cite{boardmanvogt}.  In each case the idea was to mirror the homotopy coherent algebraic structures that occurred in loop spaces,  etc. 

As an example of the problem, Segal mentions the following:  Suppose $\mathcal{C}$ is a category and that coproducts exist in $\mathcal{C}$.  How is this reflected in the nerve of $\mathcal{C}$?  It very nearly acquires a composition law, since from $X_1$ and $X_2$, one gets $X_{12} = X_1\sqcup X_2$, and two 1-simplices
$$X_1\rightarrow X_{12} \leftarrow X_2,$$
but $X_{12}$ is only determined up to isomorphism.
Let $\mathcal{C}_2$ be the category of such diagrams, i.e. in which the middle is the coproduct of the ends.  There is a functor
$$\delta_2: \mathcal{C}_2 \to \mathcal{C}\times \mathcal{C}$$and this is an equivalence of categories, but there is also a `composition law' $$m : \mathcal{C}_2 \to \mathcal{C}$$ given by picking out the coproduct. This looks fine but in fact this tentative multiplication again hits the problem of associativity.  The theory of monoidal categories  was not as developed then in 1974 as it is now, and Segal's neat solution was to side-step the issue.  He formed a category $\mathcal{C}_3$  consisting of all diagrams of form
$$
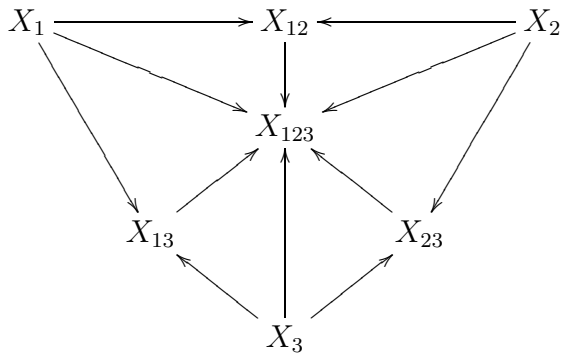
\xymatrix{X_1\ar[rr]\ar[drr]\ar[ddr]&&X_{12}\ar[d]&&X_2\ar[ll]\ar[dll]\ar[ddl]\\
&& X_{123}&&\\
&X_{13}\ar[ur]&&X_{23}\ar[ul]&\\
&&X_3\ar[ul]\ar[uu]\ar[ur]&&}$$
the notation indicating that each split line corresponds to the middle term being the coproduct of the two ends. All the usual natural isomorphisms between multiple coproducts are encoded in the one category.  There is an equivalence of categories 
$$\delta_3: \mathcal{C}_3 \to \mathcal{C}\times \mathcal{C}\times \mathcal{C}$$
sending the above diagram to $(X_1,X_2,X_3)$ and a `ternary operation' $\mathcal{C}_3 \to \mathcal{C}$ sending the diagram to $X_{123}$ compatibly, up to specifiable homotopies, with the structure outlined earlier.  The advantage is that all of this can be encoded by the nerve and thus by the classifying space structure as a $\mathbf{\Gamma}$-space. The $\mathbf{\Gamma}$-space machinery is now quite well known as it has, for instance, considerable importance in symmetric operad theory, but what are the definitions in that theory and how do they relate to our main theme.

\subsection{$\mathbf{\Gamma}$-spaces, $\mathbf{\Gamma}$-categories.}
\textbf{Definition}

(i) The category $\mathbf{\Gamma}$ is the category whose objects are all finite sets and whose morphisms from $S$ to $T$ are the maps $\theta :S \to \mathcal{P}(T)$ such that when $a\neq b\in S$ then $\theta(a)\cap \theta(b) = \emptyset$.  The composite of $\theta: S\to \mathcal{P}(T)$ and $\phi :T \to \mathcal{P}(U)$ is $\psi : S \to \mathcal{P}(U)$ where $\psi(a) = \bigcup_{b\in\theta(A)}\phi(b)$. 

(ii) A \emph{$\mathbf{\Gamma}$-space} is a functor $A: \mathbf{\Gamma}^{op}\to Top$ such that\\
(a) $A(0)$ is contractible;\\
and \\
(b) for any $n$, the map $p_n : A(n) \to\underbrace{ A(1)\times \cdots \times A(1)}_n$ induced by the maps $i_k : 1\to n$ in $\mathbf{\Gamma}$ where $i_k(1) = \{k\}\subset n$, is a homotopy equivalence.

\medskip

There is an obvious functor $\mathbf{\Delta} \to \mathbf{\Gamma}$ which takes $[m]$ to the set $m$ and $f : [m]\to [n]$ to $\theta (i) = \{j \in n  ~|~f(i-1)< j \leq f(i)\}$.  Composing a $\mathbf{\Gamma}$-space $A$ with (the opposite of) this functor gives an underlying simplicial space for the $\mathbf{\Gamma}$-space. (As we have tended to concentrate on simplicial theory rather than on the topological side, we could equally well have replaced `space' as meaning topological space by `space' as meaning simplicial set, in which case a $\mathbf{\Gamma}$-space would be a bisimplicial set with side conditions.)

Our task is not to review the contents of Segal's 1974 paper, so the next point to note is the definition of a $\mathbf{\Gamma}$-category.  This follows the same model:

A $\mathbf{\Gamma}$-category is a functor $\mathcal{C}: \mathbf{\Gamma}^{op}\to Cat$ such that (a) $\mathcal{C}(0)$ is equivalent to a one arrow category, and (b) as before except weak homotopy equivalence is replaced by equivalence of categories.

\medskip

It is not surprising that if $\mathcal{C}$ is a $\mathbf{\Gamma}$-category, applying nerve and then geometric realisation (i.e. taking its classifying space) gives a $\mathbf{\Gamma}$--space.

The most fundamental $\mathbf{\Gamma}$-category is when $Sets_{fin}$ is the category of finite sets under disjoint union, and we take an object $n$ for each natural number $n$.  The resulting $\mathbf{\Gamma}$-category is closely related to the disjoint union of the symmetric groups and thus to free symmetric monoidal closed categories, but exploring in that direction would take us too far afield.

\medskip

\textbf{Remarks}

(i) It is not obvious to start with why the structure of a $\mathbf{\Gamma}$-space is built on $\mathbf{\Gamma}$ and not just on $\mathbf{\Delta}$.  The point is, and this is important for the alternative ideas that we will be looking at later, if $A$ is a $\mathbf{\Gamma}$-space, we can form a classifying $\mathbf{\Gamma}$-space $BA$.  Any $\mathbf{\Gamma}$-space yields a simplicial space as above and hence a space $|A|$ by using the geometrical realisation (technically one uses the form of realisation that does not use the degeneracies, cf. \cite{segal} again). The classifying space $BA$ is given by assigning to a finite set $S$, the realisation of the $\mathbf{\Gamma}$-space $T\longmapsto A(S\times T)$.  This would not be possible if one considered just the underlying simplicial space, but suggests that the classifying space might be thought of as a bisimplicial space.

(ii)  It is worth recalling the structure of Lawvere's algebraic theories, \cite{lawvere}, for comparison.  One way to view these is to use some elementary ideas from topos theory.  The category of finite sets is given the structure of a basic algebraic site $\mathcal{T}_0$ by taking epimorphic families as covering sieves.  The minimal covering sieves induce colimit cones $1\sqcup 1 \cdots \sqcup 1 \to n$, so that sheaves on the algebraic site are graded sets $\{X(n) ~|~ n\geq 0\}$, endowed with bijections $X(n) \cong X(1)^n$ for $n\geq 0$.  An algebraic theory is then a `coproduct preserving' extension $\mathcal{T}_1$ of $\mathcal{T}_0$ and a $\mathcal{T}_1$-model or $\mathcal{T}_0$-algebra is a presheaf on  $\mathcal{T}_1$, which restricts to a sheaf on  $\mathcal{T}_0$.  (Some of the links between the modern theory of weak $\omega$-categories, and operads as a generalisation of Lawvere's algebraic theories are considered in Berger, \cite{berger}, and in related articles mentioned there. This is very relevant to the final theme of these notes, namely the link with algebraic models for homotopy types and the corresponding higher order categories.)  It is clear that Segal's $\mathbf{\Gamma}$-spaces are a lax or `up-to-homotopy' version of algebraic theories. We should also note out that there are links with operads of various types, but a discussion of these would take us too far afield.

\subsection{Segal categories}

The theory of $\mathcal{S}$-categories is too strict for convenience.  As we have seen the notion of Segal 1-category should be a weakened version of that of $\mathcal{S}$-category and we have discussed them informally earlier. A more formal treatment is given in several places.  The following is adapted from Toen's \cite{toen}.
\medskip

\textbf{Definitions}
\begin{itemize}
\item  A \emph{Segal 1-precategory} is specified by a functor
$$A : \mathbf{\Delta}^{op} \to \mathcal{S}$$
(i.e. a bisimplicial set) such that $A_0$ is a constant simplicial set called \emph{the set of objects of} $A$.
\item A \emph{morphism} between two Segal 1-precategories  is a natural transformation between the functors from $ \mathbf{\Delta}^{op}$ to $\mathcal{S}$.
\item A Segal 1-precategory $A$ is a \emph{Segal 1-category} if for each $[p]$, the Segal map
$$\delta[p]: A_p \to A_1\times_{A_0}A_1  \times_{A_0}\ldots \times_{A_0}A_1,$$
is a  homotopy equivalence of simplicial sets.
\item For any Segal 1-category $A$, we define its \emph{homotopy category} $Ho(A)$ as the category having $A_0$ as its set of objects and $\pi_0(A_1)$ as its set of morphisms.
\item A morphism of Segal 1-categories $f : A \to B$ is an \emph{equivalence} if it satisfies the following conditions:
\begin{enumerate}
\item For each $[p]\in \mathbf{\Delta}$, the morphism $f_p : A_p \to B_p$ is an equivalence of simplicial sets;
\item The induced functor $Ho(f) : Ho(A) \to Ho(B)$ is an equivalence of categories.
\end{enumerate}
\end{itemize}
The category of Segal 1-precategories will be denoted $1-PrCat$.  It contains the category of (small) $\mathcal{S}$-categories as we have seen. It also has a Quillen closed model category structure in which the cofibrations are the monomorphisms, and the fibrant objects are exactly the Segal 1-categories, 
moreover a morphism between Segal 1-categories is a weak equivalence in $1-PrCat$ if and only if it is an equivalence of Segal 1-categories in the above sense.  (For more on this structure, see the paper by Simpson, \cite{carlos9704006}.)

\textbf{Remark.}

Although quite a useful intermediate concept. I feel that `1-precategory' seems a slight misnomer, it is as if we called a directed graph a `precategory.  There is no algebraic structure involved in the  notion.  I have stuck with the terminaology for want of a better term.

We have seen that $\mathcal{S}$-categories and Segal categories model parts of homotopy theory well.  At the same time $G(K)$ is a $\mathcal{S}$-groupoid and we noted that these model all homotopy types. In any of these models, it is always relatively easy to extract information in low dimensions, but the level of `algebraicity' in such models is limited and many people have searched for $n$-categorical models of homotopy theory or of homotopy types that incorporate both geometric and algebraic aspects of the theory.  Based on Segal-categories and his delooping machine of \cite{segal}, one gets the Tamsamani models, \cite{tam}, for weak $n$-categories and weak $n$-groupoids for any $n$.

\section{Tamsamani weak $n$-categories} 

The problem of finding good $n$-categorical models for the geometric data encoded in a homotopy type came to the fore with Grothendieck's notes, \cite{stacks} and in particular his letter to Quillen therein.  In this he suggested that there should be models of homotopy types that (i) were less redundant in their `data storage' than, say, simplicial sets or simplicial groups, (ii) had the advantages of the mixed algebraic combinatoric power of categories (cf. the way in which the fundamental groupoid of a simplicial complex encodes the graph theoretic structure of the 1-skeleton and the algebraic nature of path concatentation), (iii) had a highly developed homotopy theory that was consistent with enough homotopy coherence and category theory to enable homotopy analogues of sheaf theoretic constructions to be made (stack theory), and (iv) for which could be developed a higher order version of the unified Galois-Poincar\'e theory that encompasses both classical Galois theory and the Galois correspondence that classifies covering spaces in terms of $\pi_1(X)$--sets, (cf. Borceux and Janeldize, \cite{BJ}). 

There are now well known problems in using the hierarchy of strict $n$-categories for such a task, as the strictness kills certain structure that is needed if one is to model spaces which have, for instance, non-trivial Whitehead products. (The crossed complex models of \cite{rb&pjh1987},  for example, are able to have a relatively simple structure because they do not try to model Whitehead products. They can thus be thought of either as incomplete models for all homotopy types or complete models for a restricted class of homotopy types, namely those with trivial Whitehead products.) From a purely categorical position, the problem arises in dimension 3, where 3-groupoids cannot capture all homotopy 2-types.  The problem is the interchange law. We saw how it corresponded to two homotopies `sliding' over each other, but that homotopy commutative square is not a trivial one and \emph{does} encode structure.   One solution is to used Gray groupoids, (see the discussion and references in \cite{K&P:gray}, for instance). Another related one is to use tri-groupoids,  the next level up from bicategories and bigroupoids. Tamsamani's approach uses Segal category ideas to encode this necessary lack of `strictness' and obtains weak $n$-categories generalising bicategories for all $n$..

\subsection{Bisimplicial models for a bicategory. }
As a starting point for our look at Tamsamani's models, consider the 2-dimensional analogue of the Grothendieck-Segal condition given in Proposition \ref{GrotSegal}.  That showed that if $A$ was a simplicial set such that for all $p\ge 2$, the Segal maps
$$\delta[p]: A_p \to A_1\times_{A_0}A_1  \times_{A_0}\ldots \times_{A_0}A_1,$$
were bijections, then $A$ was the nerve of a category. (We will just say that $A$ is a \emph{1-nerve}.)

Going up one dimension, a useful process is `categorification'.  This horrible word is in fact a good term for an important process.  A monoid is a set with some extra structure.  Categorify it and it is a small category with a single object.   A monoidal category is a category with a multiplication.  Categorify it and it is a bicategory with a single object.  In the categorification process another idea also occurs. Equality is a great idea for elements, but a poorly behaved one for monoids where isomorphism is more natural. Categorify and isomorphism is not so natural, equivalence between categories is what is needed. As one categorifies the key to success is to use equivalences except possibly in the top dimension available.

Categorifying the Grothendieck-Segal condition, let $A$ be a bisimplicial set. We can think of this as 
$$A: (\mathbf{\Delta}^{op})^2\to Sets$$
or
$$A : \mathbf{\Delta}^{op} \to \mathcal{S}.$$
Both ways are useful, but a word is needed on notation. If $A : \mathbf{\Delta}^{op} \to \mathcal{S},$ then for $[p] \in \mathbf{\Delta}$, $A_p : \mathbf{\Delta}^{op}\to Sets$ and its value on $[q]$ is $A_{p,q}$ in the other notation. Confusion can arise about which `variable' is changing and which constant, so the notation $A_{p/}$ has been introduced to indicate that  anything after the / is varying.  In dimension 2, this is not that essential but in higher dimensions it is a great help.

Assume that (i) $A_{0/}$ is a constant simplicial set, written simply $A_0$ and called the \emph{set of objects} of $A$, (ii) for each $p \in \mathbf{\Delta}$, $A_{p/} \in \mathcal{S}$ is a 1-nerve, and (iii), again for each $p \in \mathbf{\Delta}$, the morphism
$$\delta[p]: A_{p/} \to A_{1/}\times_{A_0}A_{1/}  \times_{A_0}\ldots \times_{A_0}A_{1/},$$
is an equivalence of categories.

In other words $A$ is a `Segal 2-category' or perhaps `Tamsamani-Segal weak 2-category'.  What does such a beastie look like?

For any pair $(p,q) \in \mathbb{N}^2$, we can think of any $a\in A_{p,q}$ as being a $(p,q)$-prism, $$a: \Delta[p]\times\Delta[q]\to A.$$  
The condition that $A_{0/}$ is constant, implies that $A_{0,q}$ is the same for all $q$.  For instance, a $(1,1)$-prism would naturally look like a square, but since it is constant in the second direction we get as a better picture:
\begin{center}
\includegraphics{glob1a.epsi}
\quad fig. 1
\end{center}
The iterated face maps from $A_{p/}$ to $A_0$ yield a map
$$A_{p/} \to A_0\times \ldots A_0$$
and if $(x_0, \ldots, x_p) \in A_0\times \ldots A_0$, then we may think of the inverse image of $(x_0, \ldots, x_p)$ by this map as the 1-nerve (and thus category) of those prisms having the $(x_0, \ldots, x_p)$ as their vertices in the first direction.  We will denote this by $A_{p/}(x_0, \ldots, x_p)$. In our simple case we have $p = 1$, and $A_{1/}(x_0,x_1)$ consisting of things like:
\begin{center}
\includegraphics{glob2a.epsi}
fig. 2
\end{center}
The simplicial set $A_{1/}$ is a 1-nerve, so \emph{is} a category: for each $q\geq 2$, the Segal map
$$\delta[q]: A_{1,q} \to A_{1,1}\times_{A_0}A_{1,1}  \times_{A_0}\ldots \times_{A_0}A_{1,1},$$
is a bijection.  This is inherited by the individual $A_{1/}(x_0,x_1)$s.  For instance, things in $A_{1,2}$ look like:
\begin{center}
\includegraphics{glob3a.epsi}
fig. 3
\end{center}
and as the diagram is intended to indicate, composition works perfectly.

The Segal maps
$$\delta[p]: A_{p/} \to A_{1/}\times_{A_0}A_{1/}  \times_{A_0}\ldots \times_{A_0}A_{1/},$$
are all equivalences of categories. To start to understand this we look at $p = 2$. The prism for $q=1$ looks like:
\begin{center}
\includegraphics{glob4a.epsi}
fig. 4
\end{center}\label{cushion}
We have three `lozenge' shaped vertical faces, (each with a 2-cell/arrow going from top to bottom) and the `horizontal' triangular 2-cells (no arrow shown in the diagram).

If $A_{2/}$ was \emph{isomorphic} to $A_{1/}\times_{A_0}A_{1/}$ then the front two faces, the 2-cells $a$ and $b$, would uniquely determine the third face, and thus a composite, but all we have is that $\delta[2]$ is an \emph{equivalence} of categories.  We let $\gamma[2] : A_{1/}\times_{A_0}A_{1/}\to A_{2/}$ be an inverse equivalence.  There are natural isomorphisms
$$\alpha_2: \delta[2]\gamma[2] \stackrel{\cong}{\to} Id, \mbox{\quad and \quad }\gamma[2]\delta[2] \stackrel{\cong}{\to} Id.$$
The obvious definition for the composite would be  $d_1^1\gamma[2](a,b)$, but it is better to preprocess things a bit as this does not give quite the right answer! Suppose that 
within $A_{1/}(x_0,x_1)$,  $a: f \Longrightarrow f^\prime$ and that within $A_{1/}(x_1,x_2)$,  $b: g \Longrightarrow g^\prime$. (We will use `2-cell' arrows since that is what we hope they will be later.) There are isomorphisms within $A_{1/}\times_{A_0}A_{1/}$,
$$\alpha_2(f,g) : \delta[2]\gamma[2](f,g)\Longrightarrow (f,g),$$and
$$\alpha_2(f^\prime,g^\prime) : \delta[2]\gamma[2](f^\prime,g^\prime)\Longrightarrow (f^\prime,g^\prime),$$
and hence a composite 
$$\delta[2]\gamma[2](f,g)\stackrel{\alpha_2(f,g)}{\Longrightarrow} (f,g)\stackrel{(a\times b)}{\Longrightarrow}  (f^\prime,g^\prime)\stackrel{\alpha_2(f^\prime,g^\prime)^{-1}}{\Longrightarrow}\delta[2]\gamma[2](f^\prime,g^\prime).$$
Writing $\sigma = \gamma[2](f,g)$ and similarly $\sigma^\prime = \gamma[2](f^\prime,g^\prime)$, we have two objects of the category  $A_{2/}$ and as $\delta[2]$ is an equivalence, it is fully faithful and essentially surjective.  However we have  our composite in $A_{1/}\times_{A_0}A_{1/}(\delta[2]\sigma,\delta[2]\sigma^\prime)$, so there is a unique $\varepsilon  \in A_{2/}(\sigma,\sigma^\prime)$ mapped down to it by $\delta[2]$. We set $f\#_0g = d_1^1(\sigma)$, and $f^\prime\#_0g^\prime = d_1^1(\sigma^\prime)$, giving a sensible definition of \emph{horizontal} composition of 1-arrows, and then $a\#_0 b = d^1_1(\varepsilon)$ as composite 2-cell.

 \medskip

We will not give the complete proof of the fact that any Tamsamani-Segal 2-category defines a bicategory and \emph{vice versa}.  The above just gives some of the flavour of the method. In this approach it is very easy to make a slip so let us just note that the arguments being used are more or less identical to those used to prove the result on replacement up to specified homotopy equivalence of spaces in a commutative diagram so as to get a homotopy coherent diagram. It should thus be possible to use ideas from that area of homotopy coherence, and experience of similar problems in other areas to shorten  the proof given by Tamsamani (the details are omitted from \cite{tam}), but are included in the thesis on which that article is based.

\subsection{Tamsamani-Segal weak $n$-categories}

We will not give a detailed treatment of the extension of these ideas to $n$ dimensions. It would take up too much space here. Rather we will follow the summary of that extension given by Simpson in \cite{carlos9708010}. (One word of warning, in that source the only types of $n$-categories that occur are the weak $n$-categories being developed in that setting, so the author omits the adjective `weak' consistently. He states this early on but it means that it is dangerous to `dip' into this paper although well worth reading from the start for the ideas it discusses and develops. We will call these objects `T-S weak $n$-categories' to distinguish them from the other models available.) We will launch in at `the deep end'!

\medskip

\textbf{Definition}

A \emph{T-S weak $n$-category} is a functor 
$$A : \mathbf{\Delta}^{op} \to weak(n-1)Cat,$$
where $weak(n-1)Cat$  is the category of T-S weak $(n-1)$- categories, such that
\begin{itemize}
\item $A_0$ is a set (i.e. a discrete T-S weak $(n-1)$-category);
\item  for each $p \geq 2$, the Segal map
$$\delta[p]: A_{p/} \to A_{1/}\times_{A_0}A_{1/}  \times_{A_0}\ldots \times_{A_0}A_{1/},$$
is an $(n-1)$-equivalence of weak $(n-1)$-categories.
\end{itemize}

\medskip

\textbf{Remarks}

Clearly some remarks are called for:\\
(i) It is clear that any T-S weak $n$-category is a special $n$-simplicial set and could be equally well specified  by
$$A :  (\mathbf{\Delta}^{op})^n\to Sets.$$
The notation $A_{p/}$ is the obvious extension of that which was used in case $n = 2$.  Thus the $(n-1)$-simplicial set $A_{p/} $ satisfies
$$A_{p/}([q_1], [q_2], \ldots, [q_{n-1}]) = A_{p,q_1,q_2,\ldots, q_{n-1}}.$$
A further extension of the notation is also useful.  Let $M = ([m_1], \ldots ,[m_{n-1}])$ be an object of $\mathbf{\Delta}^{n-1}$ and $[m] \in \mathbf{\Delta} $, then $(M,m)$ or more correctly $(M,[m])$, denotes an obvious object of $\mathbf{\Delta}^{n}$. We will write $A_M$ for the obvious simplicial set, in which the first $(n-1)$-variables are `clamped' at $M$. Finally if $M \in\mathbf{\Delta}^{n-h}$ and $M^\prime \in\mathbf{\Delta}^{h}$ for some $h$ then a notation $A_{(M,M^\prime)}$ will be used in the obvious way.  This is particularly useful when $M^\prime = 0_h:= (0, \ldots, 0) \in  \mathbf{\Delta}^{h}$ because of the next remark.\\
(ii) The requirement that $A_0 := A_{0/}$ be a set is to be interpreted as meaning that it is a constant $(n-1)$-simplicial set. It is not hard to show that as this will be the case for all intermediate definitions of T-S weak $k$-category, $0<k<n$, this condition implies that in the $n$-simplicial set $A$, if $p_i = 0$ for some $i$, then $A_{p_1,\ldots, p_n} = A_{p_1, \ldots, p_{i-1},0_{n-i}}$, i.e. it is independent of the values of $p_{i+1}, \ldots, p_n$. We will use the term \emph{$n$-precategory} for a $n$-simplicial set, $A :  (\mathbf{\Delta}^{op})^n\to Sets,$ which satisfies the first condition: \emph{$A_0$ is a set} and hence has the property of `constancy' mentioned above.\\
(iii) The category of T-S weak $n$-categories is defined in the obvious way, but, of course, \emph{its} definition will depend on that of weak $(n-1)$-categories and the category they form and so on.  As we have a notion of weak 2-category in this setting, that does not cause a problem. What is more of a bother is that the definition also depends recursively on a notion of equivalence of weak $(n-1)$-categories.
\medskip

\textbf{Truncation}

 A $n$-precategory will be said to be \emph{1-truncatable} if for all $M \in (\mathbf{\Delta})^{n-1}$, $A_M$ is a category (i.e. is the nerve of a category).

\medskip

If $A$ is 1-truncatable, then we can form a $(n-1)$-precategory $TA$, called the \emph{1-truncation} of $A$, which associates to $M$ the set of isomorphism classes of objects of the category, $A_M$.  If $A$ is 1-truncatable, there is a natural morphism of $(n-1)$-simplicial sets,
$$\tau(A)(M): A_{(M,0)}\to T(A)_{(M)},$$
which sends an object of the category $A_M$ to its isomorphism class.

We will say that $A$ is \emph{$k$-truncatable} for $2 \leq k \leq n$ if and only if $A$ is $(k-1)$-truncatable and $T^{k-1}A$ is 1-truncatable.
Now suppose that $A$ is $k$-truncatable, and that $1 \leq h\leq  k <n$.  Then there is a morphism of $(n-h)$-simplicial sets,
$$\tau^h(A)_M : A_{(M,0_h)}\to T^h(A)_{(M)},$$
for $M\in \mathbf{\Delta}^{n-h}$, and abbreviating $(M,0_i)$ to $M_i$, given as the composite
\begin{eqnarray*}A_{(M_h)}\stackrel{\tau(A)_{(M_{h-1})}}{\longrightarrow}T(A)_{(M_{h-1})}\stackrel{\tau(TA)_{(M_{h-2})}}{\longrightarrow}&T^2(A)_{(M_{h-2})}&\longrightarrow \ldots \\
\ldots \longrightarrow &T^{h-1}(A)_{(M_1)}&\stackrel{\tau(T^{h-1}A)_{(M})}{\longrightarrow}T(A)_{(M_{h-1})}.\end{eqnarray*}
\medskip

\textbf{$n$-equivalence}

As we have already noted, an equivalence of categories is a functor, $f : A \to B$  that is (i) fully faithful and (ii) essentially surjective.  If we want to define the notion of $n$-equivalence, one approach is to generalise these two properties. For $n=2$, $A$ and $B$ are now bicategories, `fully faithful' now  means that for each pair of objects, $x,y \in Ob(A)$. the induced functor
$$A(x,y) \to B(fx,fy)$$
is an equivalence of categories and similarly `essentially surjective' generalises to saying that every object in $B$ is equivalent to an object of the form $fx$ for some object $x$ of $A$.  With our link between T-S weak 2-categories and bicategories and also the truncation, we can see that a bicategory can be truncated twice.  If one applies $T$ once then we replace each category $A(x,y)$ by its set of isomorphism classes, thus $T(A)$ is a category and $T^2(A)$ is the set of isomorpism classes in $T(A)$. The claim is that $f : A \to B$ is essentially surjective if and only if $T^2(f)$ is surjective.  

Suppose $T^2(f)$ is surjective, and $b$ is an object of $B$, then there is an object $a$ in $A$ and an isomorphism in $T(B)$ between $b$ and $f(a)$. That interprets as saying there is a pair of 1-cells $\beta: b \to f(a)$, and $\alpha: f(a)\to b$, whose composites are isomorphic to the respective identities, i.e. there are diagrams
$$\xymatrix{&f(a)\ar[dr]^\alpha\ar@{}[d]|\Updownarrow&\\
b\ar[rr]_{id}\ar[ur]^\beta && b}\hspace{2cm}\xymatrix{&b\ar[dr]^\beta\ar@{}[d]|\Updownarrow&\\
f(a)\ar[rr]_{id}\ar[ur]^\alpha && f(a)},$$
which seems like a pretty good version of equivalent objects. The converse is similar.

Tamsamani, \cite{tam}, generalises this idea to essential surjectivity in higher dimensions, and, as Simpson notes in \cite{carlos9708010}, this can be viewed as saying that $f: A \to B$ between two T-S weak $n$-categories is a essentially surjective if $T^n(f)$ is surjective.  In fact, he points out that if $f$ is an $n$-equivalence, then $T^n(f)$ is a bijection.  

Thus we have as a definition that a morphism $f: A\to B$ of T-S weak $n$-categories is an $n$-equivalence if and only if (i) each pair of objects, $x,y \in Ob(A)$. the induced morphism
$$A_{1/}(x,y) \to B_{1/}(fx,fy)$$
is an $(n-1)$-equivalence of weak  $(n-1)$-categories, and (ii) $T^n(f)$ is surjective.

\medskip

This essentially finishes the definition of T-S weak $n$-category as we now have a working definition for all the terms involved.

It is amusing and quite useful that $T$ gives a functor from $weak~n~Cat$ to $weak~(n-1)~Cat$.

\subsection{The Poincar\'e weak $n$-groupoid of Tamsamani}

Tamsamani defines a fundamental  Poincar\'e $n$-groupoid for an arbitrary space $X$ and any $n$. There are some errors in the published version, since an attempt at an iterative definition fails since no topology has been specified on  the set of simplices involved, however this slip is rectified later on and the full approach would seem to give a fairly clear method for defining the gadget.  It has to be remembered that the result will be a weak $n$-groupoid, so in dimension $n=2$ we get a bigroupoid and not a 2-groupoid or a double groupoid with connection.  This means that the object is rather large.  Here we will attempt to describe the bigroupoid by taking apart the construction in that case.  We initially give the construction in general as it is easier to do it that way.

Let $X$ be a space and $M = (m_1, \ldots, m_n)$,  an object of $\mathbf{\Delta}^n$.  Let $\Delta ^M : = \Delta^{m_n} \times \dots \times \Delta^{m_1}$, (note the reversal of order). If $0\leq k\leq n-1$, set $M_{n-k} = (m_1, \ldots, m_{n-k-1})$ and $M^\prime_{n-k} = (m_{n-k+1}, \ldots , m_n)$  and for each $0\leq i\leq m_{n-k}$, let $v_i$ denote the $i^{th}$ vertex of $\Delta^{m_{n-k}}$.

Let
\begin{eqnarray*}X_M = \{f: \Delta^M \to X\hspace{-.3cm} &|&\hspace{-.3cm} \mbox{ for all }k, ~ 0\leq k\leq n-1,  \mbox{ for all }i,~ 0\leq i\leq m_{n-k}\\ &&\hspace{1cm}\mbox{ and for all }\mathbf{x}\in \Delta^{M_{n-k}}, \mathbf{x}^\prime \in \Delta^{M^\prime_{n-k}}, \\
&&\hspace{3cm}f(\mathbf{x},v_i,\mathbf{x}^\prime) = f_i(\mathbf{x^\prime}) \mbox{ where }f_i \in X_{M_{n-k}}\}\end{eqnarray*}
Note that if $m_{n-k} = 0$, then there is one vertex only of the corresponding $\Delta^0$ so then $X_M = X_{M_{n-k}}$.
This encodes into a subcomplex of the $n$-simplicial singular complex, the constancy rule that we need for a weak $n$-category.

\medskip

\textbf{Homotopy in $X_M$}

Let $f,g$ be two elements of $X_M$.  We will say that $f$ and $g$ are homotopic and write $f\simeq g$ or $\overline{f} = \overline{g}$ if there is a $\gamma \in X_{M,1}$ such that $\delta_0(\gamma) = f$ and $\delta(\gamma) = g$. Homotopy is an equivalence relation and we will denote by $\overline{X}_M$ the set of homotopy classes.

With the obvious identifications, we now define, $\Pi_n(X)_M := \overline{X}_M$ with the induced face and degeneracy maps.    

\begin{theorem}\cite{tam}\\
The $n$-simplicial set $\Pi_n(X)$ is a T-S weak $n$-groupoid. \hfill $\blacksquare$
\end{theorem}

\textbf{So what does it look like in dimension 2?}

The definition of $X_M$ has quite a few subcases even with $n = 2$, so we take them one at a time. We have `for all $k$, $0\leq k\leq 1$':

$k = 0$: \quad for all $i$, $0\leq i\leq m_2$, $v_i$, the $i^{th}$ vertex of $\Delta^{m_2}$, and for all $x_1 \in \Delta^{m_1}$, one has $f(v_i,x_1) = f_i(x_1)$.  However this places no restriction on $f$ since there are no variables in front of the $v_i$. It merely states that $f_i = f(v_i, \_)$.

$k = 1$: \quad for all $i$, $0\leq i\leq m_1$, $v_i$, the $i^{th}$ vertex of $\Delta^{m_1}$, and for all $x_2 \in \Delta^{m_2}$, one has $f(x_2,v_i) = f_i$, but here $f_i$ has \emph{no} variables, it is a constant value.  The picture for $M = (1,1)$ is of a square with constant values on the vertical sides as in our discussion of the bicategory model earlier in these notes. 
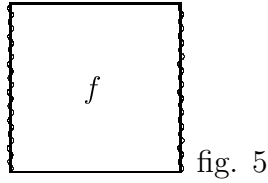
\begin{center}
\includegraphics{square1a.epsi}
fig. 5
\end{center}
The picture for an element in $X_{(2,1)}$ is similar, a vertical prism, $\Delta^2\times \Delta^1$, $\Delta^2$ as base and with `constant' vertical edges, i.e. like the shape we saw earlier in the discussion of bicategories, i.e. fig 4 on page \pageref{cushion}. The case $M=(1,2)$ is a horizontal prism with constant ends.

To study the homotopy relation, we need to look at $X_{M,1}$.  

In particular we look at $X_{1,1,1}$: there are three cases $k = 0,1,2.$ As above the case $k = 0$ imposes no condition on the singular multi-prism $f$, since the rule merely states $f(v_i,\_,\_) = f_i(\_,\_)$ in a sense just defining $f_i$.  The condition for $k=1$, however implies that for all vertices $v_i$ of $\Delta^{m_2}$,
$$f(x_3,v_i,x_1) = g_i(x_1),$$
i.e. is independent of the variable from $\Delta^{m_3}$.  Finally for $k = 2$, the restriction is that $f(x_3,x_2,v_i)$ is a constant function.  This means that a homotopy can be represented as a singular cube in $X$ in which the left and right vertical faces are constant, the top and bottom are independent of the third direction and the other two faces have no other restriction, i.e. the homotopy is precisely a homotopy relative to the boundary of the squares in dimension $(1,1)$. This makes it look as if the Tamsamani bigroupoid is essentially the same as that of Hardie, Kamps and Kieboom, (HKK), \cite{HKK:bigroupoid}. It would be interesting to see if Tamsamani's weak 3-groupoid could be adapted, as can this bigroupoid case, to allow for a notion of thinness followed by a quotient to, say, a Gray groupoid, analogously to the way in which the HKK bigroupoid leads in \cite{HKK:2-groupoid}, to a 2-groupoid, or a double groupoid with connections.

\section{Conclusion?}
 We have looked at the way in which $\mathcal{S}$-categories, their homotopy coherent form, Segal 1-categories, and, perhaps, their iterated form, the T-S weak $n$-categories, enter into the two key areas of abstract homotopy theory. Some of the sources used have been fairly recent, so there is still a lot to do.  Here is a list of some questions, some better than others:

1. What is the precise link between the Dwyer-Kan $\mathcal{S}$-groupoid and simplicial coherence?  What do homotopy coherent simplices in $G(K)$ tell one about the models? Do they lead to good descriptions of higher interchange elements analogously to the way in which maps from $S[4]$ to a $\mathcal{S}$-category produce that interchange square?  (The work of Ali Mutlu and myself on higher order Peiffer pairings in simplicial groups may be of relevance here, cf, \cite{atmp1,atmp2}.)

2. The Tamsamani method starts with a space $X$ and produces a multi-simplicial singular complex.  That method could be applied to other types of object, for example, a simplicial group, a category, and so on.  What do the corresponding Poincar\'e weak $n$-groupoids tell one? (Remember that categories model \emph{all} homotopy types.  They just do it in rather a difficult way from the point of view of calculations!)

3. Is it true that the Dwyer-Kan hammock localisation of $\mathcal{B}^\mathbb{A}$ with respect to level homotopy equivalences is `closely related to' $\mathcal{S}(S(\mathbb{A}), Ner_{h.c.}(\mathcal{B}))$ for $\mathcal{B}$ locally Kan? If so `how close'? A lot of light on this problem has been shed by Vogt in \cite{vogt92} in the topological setting. The question would seem to be of particular importance given the upsurge of interest in $A_\infty$-categories resulting from new approaches to quantum deformation.

4.  Can one construct a homotopy coherent nerve for a general Segal category (probably yes) and what are its properties?  In general what is the precise relationship between quasi categories (as a weakening of categories) and Segal categories (also a weakening of categories)? (This question is vague, of course, and would lead to many interpretations.)

5. Can one unpack a T-S weak 3-category in a sensible way? Is it sensible to try?.

6.  How powerful is the DKS coherence theorem for Segal categories? Can a clear `stand-alone' proof be given that does not depend on a lot of extra machinery? How constructive can it be made?

 7. If you take a T-S weak $n$-category as an $n$-simplicial set and extract (i) its diagonal and (ii) its Artin Mazur codiagonal, what does the structure that results look like? Is it related, perhaps just in the groupoid case, to hypercrossed complexes or hypergroupoids, \cite{duskin}.

8. To complete the `Grothendieck programme' of pursuing stacks, (i.e. to construct (and study) stacks of models for homotopy $n$-types and
to prove for the locally constant $n$-stacks, some form of Galois-Poincar\'e correspondence theorem between equivalence classes of $n$-stacks and the corresponding $n+1$-type model of the base space or topos), one needs a good theory of homotopy coherence with those models. The current attempts by Simpson, Toen and others (see papers in the references) go a long way towards such a goal, and in work by others in theoretical physics, similar approaches have been tried in very special cases (abelian models, via chain complexes etc.).  These latter approaches use a lot of the machinery sketched in these notes. It is probably fair to say that all these approaches suffer from the gulf between the technical nature of the machinery and the simple intuitions behind them. 

This last question is thus to ask is it possible to give a clear intuitive approach, to say, 2-stacks using the Segal-category machinery that does not get bogged down in a technical morass of Quillen model category theory, an enormous amount of weak infinity category theory or similar machinery.  This is not to say that approaches using those ideas are not providing a necessary step on the road to understanding the Grothendieck problem, but to ask for a simple approach that will aid the geometric intuition.

T.Porter,\\
Mathematics Department,\\School of Informatics,\\University of Wales Bangor,\\
Bangor,\\Gwynedd, LL57 1UT,\\ United Kingdom.
\end{document}